\documentclass[a4paper,11pt]{article}
\usepackage{amssymb,amsmath,amscd,amsfonts,amsthm,latexsym}
\usepackage{color,graphpap,fancybox,framed,enumerate}
\usepackage[dvips]{epsfig}
\usepackage[english,activeacute]{babel}
\usepackage{hyperref}

\usepackage[right=3cm,left=3cm,top=3cm,bottom=3cm]{geometry}

\usepackage{xcolor}

\newtheorem{theorem}{Theorem}[section]
\newtheorem{prop}[theorem]{Proposition}
\newtheorem{lemma}[theorem]{Lemma}
\newtheorem{defi}[theorem]{Definition}

\newenvironment{demo}{ \noindent \emph{\textbf{Proof:}}}{\hfill$\square$\\}

\newcommand{\RR}{\mathbb{R}}
\newcommand{\NN}{\mathbb{N}}
\newcommand{\CC}{\mathbb{C}}
\newcommand{\TT}{\mathbb{T}}
\renewcommand{\SS}{\mathbb{S}}

\newcommand{\Cc}{\mathcal{C}}

\newcommand{\Ac}{\mathcal{A}}
\newcommand{\Bc}{\mathcal{B}}
\newcommand{\Fc}{\mathcal{F}}
\newcommand{\Lc}{\mathcal{L}}
\newcommand{\Pc}{\mathcal{P}}
\newcommand{\Qc}{\mathcal{Q}}
\newcommand{\Oc}{\mathcal{O}}

\newcommand{\grad}{\nabla}
\renewcommand{\div}{\operatorname{div}}

\newcommand{\no}{n$^{\text{o}}$}

\newcommand{\vect}[2]{\left(\begin{array}{c} #1 \\ #2 \end{array}\right)}
\newcommand{\trans}{{^{\intercal}}}
\newcommand{\pc}{ \usefont{T1}{cmtl}{m}{n} \selectfont}

\newdimen\texpscorrection
\newdimen\figcenter
\def\figurewithtex #1 #2 #3 #4 #5\cr{\null
  {\goodbreak\figcenter=\hsize\relax
  \advance\figcenter by -#4truecm
  \divide\figcenter by 2
  \begin{figure}[hbt]
  \vskip #3truecm\noindent\hskip\figcenter
  \includegraphics{#1}{\hskip\texpscorrection\input #2 }
  \vskip 0.8truecm{\baselineskip=0.8\baselineskip
  \noindent \vbox{\noindent {\footnotesize #5}}\par}
  \end{figure}}}
\def\point#1 #2 #3 {\rlap{\kern #1 truecm
\raise #2 truecm \hbox{#3}}}

\numberwithin{equation}{section}

\begin{document}

\title{\bf Decay of semilinear damped wave equations: cases without geometric 
control condition}

\author{ Romain \textsc{Joly}\footnote{Universit\'e Grenoble Alpes, CNRS, 
Institut Fourier, F-38000 Grenoble, France, email: {\pc 
romain.joly@univ-grenoble-alpes.fr}} {~\&~}
Camille \textsc{Laurent}\footnote{CNRS, UMR 7598, Laboratoire Jacques-Louis
Lions, F-75005, Paris, France  }
\footnote{UPMC Univ Paris 06, UMR 7598, Laboratoire Jacques-Louis Lions,
F-75005, Paris, France, email: {\pc laurent@ann.jussieu.fr }}}
\date{}

\maketitle

\begin{abstract}
We consider the semilinear damped wave equation
$$\partial_{tt}^2 u(x,t) + \gamma(x)\partial_t u(x,t) = \Delta u(x,t) - \alpha 
u(x,t) - f(x,u(x,t))~.$$
In this article, we obtain the first results concerning the 
stabilization of this semilinear equation in cases where $\gamma$ does not 
satisfy the geometric control condition. When some of the geodesic rays are 
trapped, the stabilization of the linear semigroup is semi-uniform in the 
sense that $\|e^{At}A^{-1}\| \leq h(t)$ for some function $h$ with 
$h(t)\rightarrow 0$ when $t\rightarrow +\infty$. We provide general tools to 
deal with the semilinear stabilization problem in the case where $h(t)$ has a 
sufficiently fast decay.\\[2mm]
\end{abstract}

\tableofcontents

\section{Introduction}

We consider the semilinear damped wave equation 
\begin{equation}\label{eq}
\left\{\begin{array}{ll}
\partial_{tt}^2 u(x,t) + \gamma(x)\partial_t u(x,t) = \Delta u(x,t) - \alpha 
u(x,t) - 
f(x,u(x,t))~~&(x,t)\in\Omega\times (0,+\infty)\\
u_{|\partial\Omega}(x,t)=0&(x,t)\in\partial\Omega\times (0,+\infty)\\
(u(\cdot,t=0),\partial_t u(\cdot,t=0))=U_0=(u_0,u_1)\in H^1_0(\Omega)\times 
L^2(\Omega)
\end{array}\right.
\end{equation}
in the following general framework:
\begin{enumerate}[(i)]
\item the domain $\Omega$ is a two-dimensional smooth 
compact and connected manifold with or without smooth boundary. If $\Omega$ is 
not flat, $\Delta$ has to be taken as Beltrami Laplacian operator.
\item the constant $\alpha\geq 0$ is a non-negative constant. We 
require that $\alpha>0$ in the case without boundary to ensure that 
$\Delta-\alpha Id$ is a negative definite self-adjoint operator.  
\item the damping $\gamma\in L^\infty(\Omega,\RR_+)$ is a bounded 
function with non-negative values. Since we want to consider a damped equation, 
we will assume that $\gamma$ does not vanish everywhere.
\item the non-linearity $f\in\Cc^1(\overline\Omega\times\RR,\RR)$ is 
of polynomial type in the sense that there exists a constant $C$ and a power 
$p\geq 1$ such that for all $(x,u)\in  \overline\Omega\times\RR$,
\begin{equation}\label{hyp-f-estim}
|f(x,u)|+|\grad_xf(x,u)|\leq C(1+|u|)^p~~\text{ and }~~|f'_u(x,u)|\leq 
C(1+|u|)^{p-1}~.
\end{equation}
Moreover, in most of this paper, we will be interested in the stabilization 
problem and we will also assume that 
\begin{equation}\label{hyp-f-sign}
\forall (x,u)\in  \overline\Omega\times\RR~,~~f(x,u)u\geq 0~.
\end{equation}
\end{enumerate}
We introduce the space $X=H^1_0(\Omega)\times L^2(\Omega)$ and the operator $A$ 
defined by 
$$D(A)=(H^2(\Omega)\cap H^1_0(\Omega))\times H^1_0(\Omega)~~~~~
A = \left(\begin{array}{cc} 0 & Id \\ \Delta-\alpha Id & -\gamma(x) 
\end{array}\right)~.$$
In this paper, we are interested in the cases where the linear semigroup 
$e^{At}$ has no uniform decay, that is that $\|e^{At}\|_{\Lc(X)}$ does no 
converge to zero. We only assume a semi-uniform decay, but sufficiently fast in 
the following sense.
\begin{enumerate}[(i)]
\setcounter{enumi}{4}
\item There exist a function $h(t)$ such that 
\begin{equation}\label{hyp-decay-1}
\forall U_0\in D(A)~,~\|e^{At}U_0\|_X \, \leq \, h(t)\,\|U_0\|_{D(A)}
\end{equation}
and there is $\sigma_h\in(0,1]$ such that
\begin{equation}\label{hyp-decay-2}
\lim_{t\rightarrow \infty}h(t)=0~~~\text{ and }~
\forall \sigma\in[0,\sigma_h]~,~~\int_0^\infty 
h(t)^{1-\sigma}\,{\rm d}t~<~\infty~. 
\end{equation}
\end{enumerate}
Condition \eqref{hyp-decay-2} requires a decay rate fast enough to be 
integrable. Roughly speaking, this article shows that this condition, together 
with a suitable unique continuation property, are sufficient to obtain a 
stabilization of the semilinear equation. The relevant unique continuation 
property is explained in Proposition \ref{prop-UCP-implique-th} below. We 
present two general results where it can be obtained.

\begin{figure}[tp]              
\ovalbox{
$~$\\
\begin{minipage}{14.7cm}
\begin{minipage}{5cm}
\epsfig{width=5cm,file=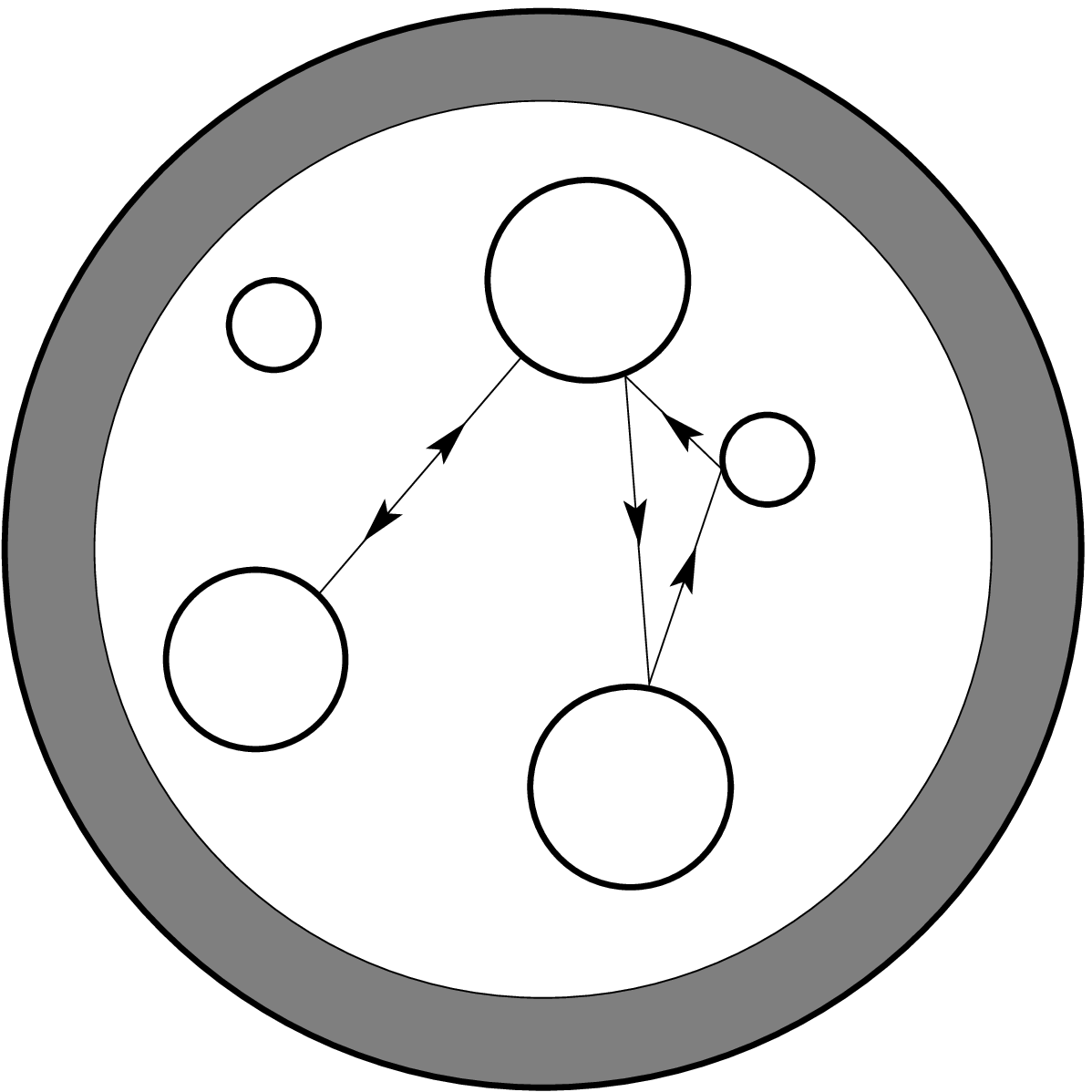}
\end{minipage}
~
\begin{minipage}{9.5cm}
{\bf\noindent Disk with holes}\\
We set $\Omega$ to be a convex flat surface with a damping 
$\gamma$ efficient near the boundary. Typically, we may take the flat disk 
$B(0,1)$ of $\RR^2$ and assume that there exist $r\in(0,1)$ and 
$\underline\gamma>0$ such that $\gamma(x)\geq \underline\gamma>0$ for $|x|>r$.
Inside the interior zone without damping, assume that at least two holes exist; 
to simplify, we also assume that these holes are disks and are small in a sense 
specified later. Notice that there exist some periodic geodesics which do not 
meet the support of the damping. This example has been studied in 
\cite{Burq-memoire, Burq-Zworski}.
\end{minipage}
\end{minipage}$~~$}

\vspace{1mm}

\ovalbox{
$~$\\[-1mm]
\begin{minipage}{14.7cm}
\vspace{1mm}
\begin{minipage}{6.5cm}
{\bf\noindent The peanut of rotation}\\
We consider a compact two-dimensional manifold without boundary. We assume that 
the damping $\gamma$ is effective, that is uniformly positive, everywhere 
except in the central part of the manifold. This part is a manifold of negative curvature and
\end{minipage}
~~
\begin{minipage}{7.6cm}
\epsfig{width=7.6cm,file=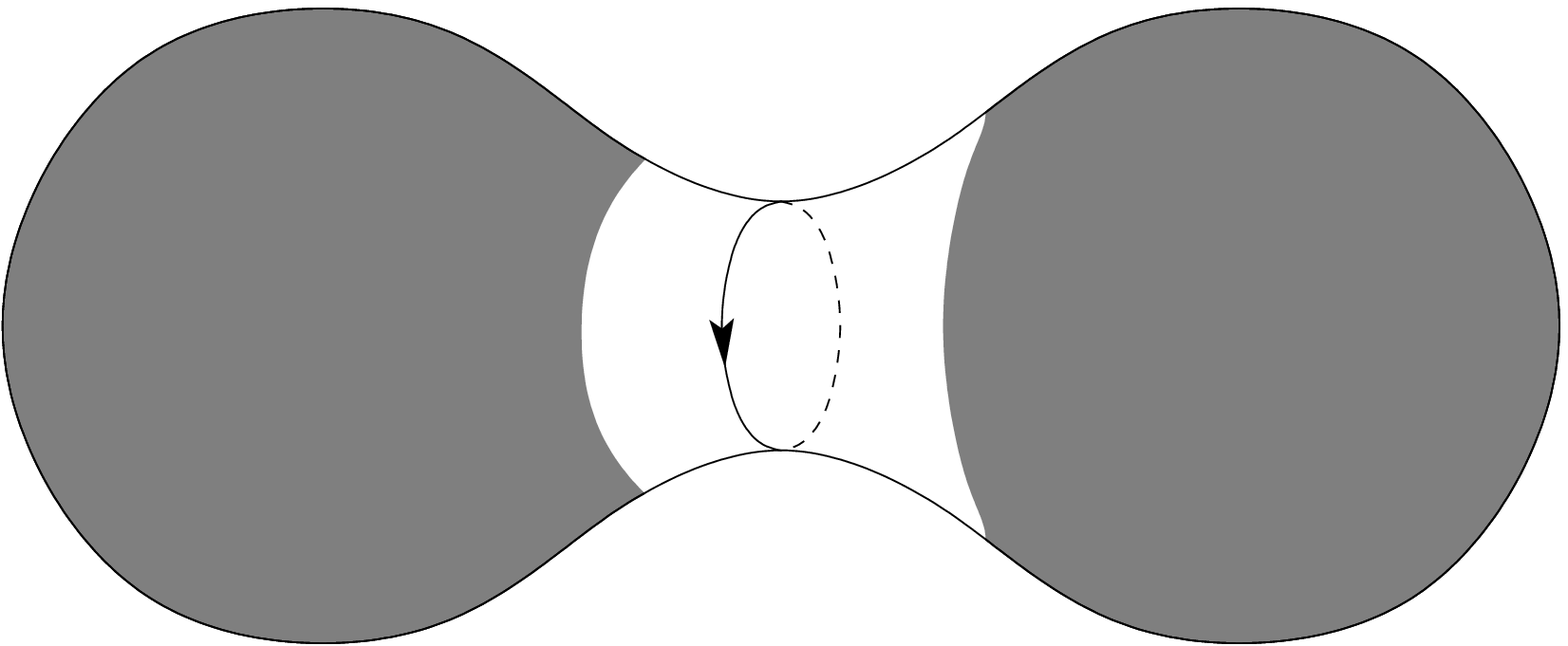}
\end{minipage}\\[1mm]
invariant by rotation along the $y-$axis. 
More precisely, let us set this part to be equivalent to the cylinder endowed 
with the metric $g(y,\theta)={\rm d}y^2+\cosh^2(y){\rm d}\theta^2$.
This central part admits a unique (up to change of 
orientation) periodic geodesic  which is 
unstable; any other geodesic meets the support of the damping. 
This example has been studied for example in \cite{CSVW}, \cite{Schenck} and 
the references therein.\\[-4mm]
\end{minipage}$~~$}

\vspace{1mm}

\ovalbox{
$~$\\
~\begin{minipage}{14.7cm}
\begin{minipage}{4.8cm}
\resizebox{4.2cm}{!}{\input{open-book.pstex_t}}
\end{minipage}
~
\begin{minipage}{9.5cm}
$~$\\
{\bf\noindent The open book}\\
We consider the torus $\TT^2$ with flat geometry. The damping $\gamma$ is 
assumed to depend only on the first coordinate and to be of the type 
$\gamma(x)=|x_1|^\beta$ with $\beta>0$ to be chosen small enough. In this case, 
there is a unique (up to change of orientation) geodesic which does not meet 
the support of the damping. This example has been studied in 
\cite{Leautaud-Lerner}.\\
\end{minipage}
\end{minipage}$~~$}

\vspace{1mm}

\ovalbox{
$~$\\
~\begin{minipage}{14.7cm}
\begin{minipage}{7.5cm}
\vspace{2mm}

{\bf\noindent Hyperbolic surfaces}\\
We consider a compact connected hyperbolic surface with constant negative curvature -1 
(for example a surface of genus 2 cut out from Poincar\'e disk). 
The damping is any non zero function $\gamma(x)\geq 0$. This example has been studied in 
\cite{Lin} following the fractal uncertainty principle of \cite{BourgainDyatlov}. 
We also refer to other results in any dimension with pressure conditions \cite{Schenckpressure} following 
ideas of \cite{Anant}.\\
$~$
\end{minipage}
~~
\begin{minipage}{7cm}
\epsfig{width=7cm,file=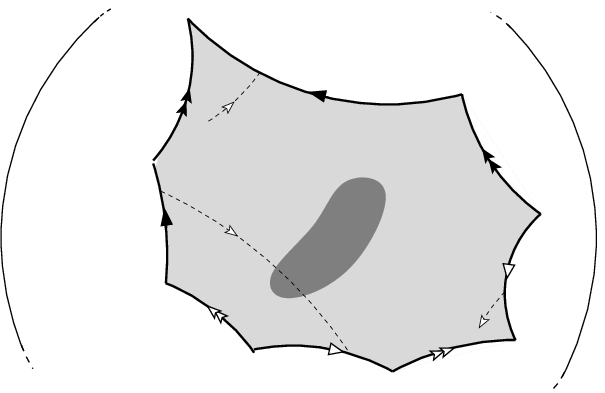}
\end{minipage}
\end{minipage}$~~$}

\caption{\it the main applications of Theorems \ref{th1} and \ref{th2} 
presented in this paper. The gray parts show the localization of the damping 
(white=no damping). The more geometrically constrained Theorem \ref{th2} apply 
to the ``disk with two holes'', the ``peanut of rotation'' and the ``open 
book''.  \label{fig-1}}
\end{figure}

Our first result concerns analytic nonlinearities and smooth dampings.
\begin{theorem}\label{th1}
Consider the damped wave equation \eqref{eq} in the framework of Assumptions 
{\normalfont (i)-(v)}. Assume in addition that:
\begin{itemize}
\item[a)] the function $(x,u)\mapsto f(x,u)$ is smooth and analytic 
with respect to $u$. 
\item[b)] the damping $\gamma$ is of class $\Cc^1$
or at least that there exists $\tilde\gamma\in\Cc^1(\overline\Omega,\RR_+)$ such that
{\normalfont (v)} holds with $\gamma$ replaced by $\tilde\gamma$ and such that the support 
of $\tilde\gamma$ is contained in the support of $\gamma$.
\item[c)] the power $p$ of $f$ in \eqref{hyp-f-estim} and the decay rate $h(t)$ 
of the semigroup in \eqref{hyp-decay-1} satisfy $h(t)=\Oc(t^{-\beta})$ with 
$\beta>2p$.
\end{itemize}
Then, any solution $u$ of \eqref{eq} satisfies 
$$\|(u,\partial_t u)(t)\|_{H^1_0\times L^2}  ~\xrightarrow[~~t\longrightarrow 
+\infty~~]{}~0~.$$
Moreover, for any $R$ and $\sigma>0$, there exists 
$h_{R,\sigma}(t)$ which goes to zero when $t$ goes to $+\infty$ such that 
the following stabilization hold. For any $U_0\in 
H^{1+\sigma}_0(\Omega)\times H^\sigma(\Omega)$, if $u$ is the 
solution of \eqref{eq}, then 
$$ \|(u_0,u_1)\|_{H^{1+\sigma}\times H^\sigma}\leq 
R~~\Longrightarrow ~~ \|(u,\partial_t u)(t)\|_{H^1_0\times L^2} \leq 
h_{R,\sigma}(t) ~\xrightarrow[~~t\longrightarrow +\infty~~]{}~0~.$$
\end{theorem}
Our assumptions (v) and c) on the decay of the linear semigroup may seem 
strong. They are satisfied in the cases where the set of trapped geodesics, the 
ones which do not meet the support of the damping, is small and 
hyperbolic in some sense. Several geometries satisfying (v) and c) have been 
studied in the literature, see the concrete examples of Figure \ref{fig-1} and 
the references 
therein. Notice in particular that the example of domain with holes is 
particularly relevant for applications where we want to stabilize a nonlinear 
material with holes by adding a damping or a control in the external part.  
There is a huge literature about the damped wave equation and the
purpose of the examples presented here is mainly to illustrate our theorem to 
non specialists. Moreover, 
the subject is growing fastly, giving more and more examples of geometries where 
we understand the effect of the damping and where we may be able to apply
our results. We do not pretend to exhaustivity and 
refer to the bibliography of the more recent \cite{Lin} for instance.

In some cases, the unique continuation property required in Proposition 
\ref{prop-UCP-implique-th} can be obtained without considering analytic 
nonlinearities or conditions on the growth of $f$ as Hypothesis c) of 
Theorem \ref{th1}. Instead, we require a particular geometry, 
which will be introduced more precisely in Section \ref{section-UCP}.
\begin{theorem}\label{th2}
Consider the damped wave equation \eqref{eq} in the framework of Assumptions 
{\normalfont (i)-(v)}. Assume in addition that:
\begin{itemize}
\item[a)] the function $(x,u)\mapsto f(x,u)$ is of class 
$\Cc^1(\overline{\Omega}\times \RR,\RR)$.
\item[b)] there exists a pseudo-convex foliation of $\Omega$ in the sense of 
Definitions \ref{defi-foliation1} or \ref{defi-foliation2}.
\end{itemize}
Then, the conclusions of Theorem \ref{th1} hold.
\end{theorem}
This result can be applied in several situations of Figure \ref{fig-1}: 
the ``disk 
with two holes'', the ``peanut of rotation'' and the ``open book''. In these 
cases, the stabilization holds for any natural nonlinearity.

We expect that the decay rate $h_{R,\sigma}(t)$ is related to the linear decay rate 
$h(t)$ of Assumption (v). We are able to obtain this link for the typical 
decays of the examples of Figure \ref{fig-1}.
\begin{prop}\label{prop1}
Consider a situation where the stabilization stated in Theorems \ref{th1} or 
\ref{th2} holds. Then, 
\begin{itemize}
\item[$\bullet$] if the decay rate of Assumption (v) satisfies 
$h(t)=\Oc({t^{-\alpha}})$ with $\alpha>1$, then the nonlinear 
equation admits a decay of the type $h_{R,\sigma}(t)= \Oc (t^{-\sigma\alpha})$.
\item[$\bullet$] if the decay rate of Assumption (v) satisfies 
$h(t)=\Oc (e^{-a {t^{1/\beta}}})$ with $a>0$ and $\beta>0$, then the nonlinear 
equation admits a decay of the type $h_{R,\sigma}(t)= \Oc (e^{-b\sigma 
t^{1/(\beta+1)}})$ for some $b>0$.  
\end{itemize}
\end{prop}
Notice that this result is purely local in the sense that the decay rate is 
obtained when the solution is close enough to $0$. Our proofs do not provide an 
explicit estimate of the time needed to enter this small neighborhood of $0$.
Also notice that the loss in the power of the second case of Proposition 
\ref{prop1} is due to an abstract setting: in the concrete examples, we may 
avoid this loss, see the remark below Lemma \ref{lemma-appli2} and the concrete 
applications to the examples of Figure \ref{fig-1}.

\vspace{3mm}

To our knowledge, Theorems \ref{th1} and \ref{th2} are the first stabilization 
results for the semilinear damped wave equation when the geometric control 
condition fails. This famous 
geometric control condition has been introduced in the works of Bardos, Lebeau, 
Rauch and Taylor (see \cite{BLR}) and roughly requires that any geodesic of the 
manifold $\Omega$ meets the support of the damping $\gamma$. This condition 
implies that the linear semigroup of the damped wave equation satisfies a 
uniform decay $\|e^{At}\|_{\Lc(H^1_0\times L^2)}\leq Me^{-\lambda t}$. In this 
context, the stabilization of the semilinear damped wave equation has been 
studied since a long time, see for example \cite{Haraux, Zuazua, Dehman, DLZ,RC}.
Under this condition and for large data, the proof often divides into a part 
dealing with high frequencies with linear arguments and another one dealing
with low frequencies that often requires 
a unique continuation argument. The high frequency problem was solved by Dehman 
\cite{Dehman} with important extension by Dehman-Lebeau-Zuazua \cite{DLZ} using 
microlocal defect measure. Yet, the unique continuation was proved by classical 
Carleman estimates (see Section \ref{subsectUCPwithout} below) which restricted 
the generality of the geometry. Using techniques from dynamical systems applied 
to PDEs, the authors of the present article proved in \cite{RC} a general stabilization 
result under Geometric Control Condition, at the cost of an assumptions of analyticity 
of the nonlinearity. Theorem \ref{th1} is in the same spirit as \cite{RC} and intends to prove 
that related techniques can be extended to a weaker damping. 
Theorem \ref{th2} is more in the spirit of the other references, taking advantage of 
particular geometries, but avoiding analyticity.

Notice that the ``disk with one hole'' satisfies this geometric control 
condition and thus it is not considered in this paper.

In the cases where the geometric control condition fails, the decay of the 
linear semigroup is not uniform. At least, if $\gamma$ does not vanish 
everywhere, it is proved in \cite{Dafermos} (see also \cite{Haraux}) that the 
trajectories of the linear semigroup goes to zero (see Theorem 
\ref{th-stab-linear} below). In fact, the decay can be estimated with a loss 
of derivative as
\begin{equation}\label{eq-weak-decay}
\|e^{At}U\|_{H^1\times L^2}\,\leq\, h(t) 
\|U\|_{H^2\times H^1}~~~~~~\text{ with }~~h(t)\xrightarrow[~t\rightarrow 
+\infty~]{} 0~.
\end{equation}
In the general case, as soon as $\gamma\not\equiv 0$, the decay rate can be 
taken as $h(t)=\Oc(\ln(\ln(t))/\ln t)$ as shown in 
\cite{Lebeau,Lebeau-Robbiano}. In some particular situations, $\gamma$ misses 
the geometric control condition but very closely: typically there is only one 
(up to symmetries) geodesic which does not meet the support of the damping and 
this geodesic is unstable. In this case, we may hope a better decay than the 
$\Oc(\ln(\ln(t))/\ln t)$ one, see for example \cite{CSVW,Leautaud-Lerner} and 
the other references of Figure \ref{fig-1}. 

To our knowledge, until now, there was no result concerning the semilinear 
damped wave equation \eqref{eq} when the geometric control condition fails. 
Thus Theorems \ref{th1} and \ref{th2} provides the first examples of 
semi-uniform stabilization for the semilinear damped wave equation. Notice that 
our results deeply rely on the fact that the decay rate of \eqref{eq-weak-decay} 
is integrable. Typically, for the situations of Figure \ref{fig-1}, it is of 
the type $h(t)=\Oc(e^{-\lambda t^\alpha})$ or $h(t)=\Oc(1/t^\beta)$ with 
sufficiently large $\beta>0$.

\vspace{3mm}

Theorems \ref{th1} and \ref{th2} concern the stabilization of the solutions of 
\eqref{eq} in the sense that their $H^1\times L^2$-norm goes to zero. Notice 
that, since 
the energy of the damped wave equation is non-increasing (see Section 
\ref{section-basic}), we knew that this $H^1\times L^2$-norm is at least 
bounded. Such a uniform bound is not clear a priori for the $H^2\times 
H^1$-norm. However, basic arguments provide this bound as a 
corollary of Theorems \ref{th1} and \ref{th2} if the decay is fast enough, 
which is the case of the ``disk with holes'', the ``peanut of rotation'' and 
the ``hyperbolic surfaces'' of Figure \ref{fig-1}.
\begin{theorem}\label{th-norm-H2}
Consider the damped wave equation \eqref{eq} in the framework of Theorems 
\ref{th1} or \ref{th2}. Assume that for all $R>0$, the decay rate $h_{R,1}(t)$ 
is faster than polynomial, i.e. $h_{R,1}(t)=o(t^{-k})$ for any $k\in\NN$. 
Also assume that $\gamma$ is of class $\Cc^1$ and $f$ is of class 
$\Cc^2(\overline\Omega\times \RR,\RR)$.
Then the $H^2(\Omega)\times H^1(\Omega)$-norm of the solutions are bounded in 
the following sense. For any $R>0$, there exists $C(R)>0$ such that, for any 
$U_0\in (H^{2}(\Omega)\cap H^1_0(\Omega))\times H^1_0(\Omega)$ such that 
$\|U_0\|_{H^2\times H^1} \leq R$, the solution $u$ of \eqref{eq} satisfies 
$$ \sup_{t\geq 0} \|(u,\partial_t u)(t)\|_{H^2\times H^1} \leq C(R)~.$$
\end{theorem}
Note that, in the case without damping, this result is sometimes expected to be 
false. It is related to the weak turbulence, described as a transport from low 
frequencies to high frequencies.

\vspace{3mm}

The main purpose of this paper is to obtain new examples of stabilization 
for the semilinear damped wave equation and to introduce the 
corresponding methods and tools. We do not pretend to be exhaustive and the 
method may be easily used to obtain further or more precise results. For 
example:
\begin{itemize}
\item the boundary condition may be modified, typically in the case of the disk 
with holes, Neumann boundary condition may be chosen at the exterior boundary.
\item for simplicity, the examples of Figure \ref{fig-1} and the main results 
of this article concern two-dimensional manifolds. However, the arguments of 
this paper can be used to deal with higher-dimensional manifolds. There are 
some technical complications, mainly due to the Sobolev embeddings. For 
example, in dimension $d=3$, the degree of $f$ in \eqref{hyp-f-estim} 
should satisfy $p<3$ ($p<5$ if we use Strichartz estimates as done in \cite{RC} 
using \cite{DLZ}) and the order $\beta$ of the vanishing of $\gamma$ in the 
example of the open book should not be too large. To simplify, we choose to 
state our results in dimension $d=2$. However, several intermediate results in 
this article are stated for dimensions $d=2$ or $d=3$. 
\item It is also possible to combine the strategy of this paper with other 
tricks and technical arguments. 
For example, we may consider unbounded manifolds or manifolds of dimension 
$d=3$ with nonlinearity of degree $p\in [3,5)$, which are supercritical in the 
Sobolev sense. This would requires to use Strichartz estimates in addition to 
Sobolev embeddings as done in \cite{DLZ} or \cite{RC}.
\item assume that we replace the sign condition \eqref{hyp-f-sign} by an 
asymptotic sign condition 
$$\exists R>0~,~~\forall (x,u)\in  
\overline\Omega\times\RR~,~~|u|\geq R ~\Rightarrow ~ f(x,u)u\geq 0~. $$
Then they may exist several equilibrium points and the stabilization to zero 
cannot be expected. However, the arguments of this paper show that the energy 
$E$ introduced in Section \ref{section-basic} is a strict Lyapounov functional 
and that any solution converges to the set of equilibrium points. We can also 
show the existence of a {\it weak compact attractor} in the sense that there is 
an invariant compact set $\Ac\subset H^1_0(\Omega)\times L^2(\Omega)$, which 
consists of all the bounded trajectories and such that any regular set 
$\Bc$ bounded in $(H^2(\Omega)\cap H^1_0(\Omega))\times H^1_0(\Omega)$ is 
attracted by $\Ac$ in the topology of $H^1_0(\Omega)\times L^2(\Omega)$.
Notice that this concept of weak attractor is the one of Babin and Vishik in 
\cite{Babin-Vishik}. At this time, the asymptotic compactness property of the 
semilinear damped wave equation was not discovered and people thought that a 
strong attractor (attracting bounded sets of $H^1_0(\Omega)\times L^2(\Omega)$) 
was impossible due to the lack of regularization property for the damped wave 
equation. Few years later, Hale \cite{Hale-book} and Haraux \cite{Haraux} 
obtain this asymptotic compactness property and the existence of a strong 
attractor. Thus this notion of weak attractor has been forgotten. It is 
noteworthy that it appears again here. Notice that we cannot hope a better 
attraction property since even in the linear case, $\{0\}$ is not an attractor 
in the strong sense.
\end{itemize}

\vspace{3mm}

The organization of this paper follows the proof of stabilization of the 
examples of Figure \ref{fig-1}. We add step by step the techniques required to 
deal with our guiding examples, from the simplest to the most complicated one. 

Sections \ref{section-basic} and \ref{section-asymp-compact} contain the basic 
notations and properties. The asymptotic compactness of the semilinear dynamics 
is proved and the problem is reduced to a unique continuation property. 
In Section \ref{section-decay}, we show the estimations of Proposition 
\ref{prop1}. Section \ref{section-appli1} then proves the nonlinear 
stabilization in the ``open book'' case, where the 
unique continuation property is trivial. Section \ref{section-UCP} stated 
several unique continuation results, enabling to prove Theorem \ref{th2}. We 
obtain as a consequence the stabilization in the case of 
the ``peanut of rotation'' in Section \ref{section-appli2}. Section 
\ref{section-decay-disk} studies the linear semigroup for the case of the 
``disk with holes'' before we apply Theorem \ref{th2} in the case of the 
``disk with two holes'' in Section \ref{section-appli3}. 
Theorem \ref{th1} is proved in Section \ref{section-reg-analytic} by showing an 
asymptotic analytic regularization. It is applied to the ``disk with three or 
more holes'' and hyperbolic surfaces, assuming $f$ analytic in $u$, in 
Sections \ref{section-appli4} and \ref{section-hyperbo}.
In Section \ref{section-th-norm-H2}, we show how to obtain Theorem 
\ref{th-norm-H2} as a corollary of Theorems \ref{th1} or \ref{th2}.
This article finishes with three appendices on the links between the decay of 
the semigroup $e^{At}$ and the resolvent $(A-i\mu Id)^{-1}$.

\vspace{8mm}

{\noindent \bf Acknowledgements: }The authors deeply thank Matthieu L\'eautaud 
for his contributions to the appendices. They are also grateful 
to Nicolas Burq for several discussions and the suggestion of Theorem 
\ref{th-norm-H2}. Part of this work has been made in 
the fruitful atmosphere of the {\it Science Center of Benasque Pedro Pascual} 
and has been supported by the project {\it ISDEEC} ANR-16-CE40-0013.

\section{Notations and basic facts}\label{section-basic}
We use the notations of Equation \eqref{eq}, of Assumptions (i)-(v) and of the 
introduction. In particular, we recall that 
$X=H^1_0(\Omega)\times L^2(\Omega)$ and 
$$D(A)=(H^2(\Omega)\cap H^1_0(\Omega))\times H^1_0(\Omega)~~~~~
A = \left(\begin{array}{cc} 0 & Id \\ \Delta-\alpha Id & -\gamma(x) 
\end{array}\right)~.$$
The operator $A$ is the classical linear damped wave operator corresponding to 
the linear part of \eqref{eq}. Due to Lumer-Phillips theorem, we know that this 
operator generates a linear semigroup of contractions $e^{At}$ on $X$ and on 
$D(A)$ and that
$$\forall t\geq 0~,~~ \|e^{At}\|_{\Lc(X)}\leq 1~\text{ and 
}\|e^{At}\|_{\Lc(D(A))}\leq 1~.$$
Notice that the second estimate is a direct consequence of the commutation of 
$A$ and $e^{At}$ and does not require any regularity on $\gamma$.

For any $\sigma\in [0,1]$, we set 
$$X^\sigma=(H^{1+\sigma}(\Omega)\cap H^1_0(\Omega))\times 
H^\sigma_0(\Omega)~.$$
Thus $X^0=X$ and $X^1=D(A)$ and $X^\sigma$ is an interpolation space 
between $X^0$ and $X^1$. In particular, by interpolation, $e^{At}$ is defined 
in $X^\sigma$ and we have
\begin{equation}\label{eq-bound-semigroup}
\forall \sigma\in [0,1]~,~~\forall t\geq 0~,~~ 
\|e^{At}\|_{\Lc(X^\sigma)}\leq 1~.
\end{equation}

We set $F\in\Cc^0(X)$ to be the function
\begin{equation}\label{def-F}
F~:~U=\vect uv \in X~\longmapsto ~ \vect 0 {-f(\cdot,u)} \in X~. 
\end{equation}
Notice that, if $\Omega$ is two-dimensional,  $H^1_0(\Omega) \hookrightarrow 
L^{p}(\Omega)$ for any $p\in[1,+\infty)$ and if $\Omega$ is three-dimensional 
$H^1_0(\Omega) \hookrightarrow L^6(\Omega)$. Thus, $f(u)$ is 
well defined in $L^2(\Omega)$ due to Assumption \eqref{hyp-f-estim} if 
dim$(\Omega)=2$ or if dim$(\Omega)=3$ and $p\leq 3$. Moreover, for any $R$, $u$ 
and $v$ with $\|u\|_{H^1}\leq R$ and $\|v\|_{H^1}\leq R$, we have
$$\|f(\cdot,u)-f(\cdot,v)\|_{L^2}=\left\|(u-v)\int_0^1 
f'_u(\cdot,u+s(u-v))ds\right\|_{L^2} 
\leq C(R)\|u-v\|_{H^1}$$
and so $F$ is lipschitzian on the bounded sets of $X$. As a consequence, 
the damped wave equation \eqref{eq} is well posed in $X$ and admits local 
solutions if dim$(\Omega)=2$ or if dim$(\Omega)=3$ and $p\leq 3$.

With the above notation, our main equation writes
\begin{equation}\label{eq2}
\partial_t U=AU + F(U)~~~~U(t=0)=U_0\in X~.
\end{equation}
In particular, Duhamel's formula yields
$$U(t)=e^{At}U_0 + \int_0^t e^{A(t-s)}F(U(s))\,{\rm d}s~.$$

We introduce the potential
$$V(x,u)=\int_0^u f(x,s)ds~.$$
Due to \eqref{hyp-f-estim} and the above arguments, $V(\cdot,u)$ defines a 
Lipschitz function from the bounded sets of $H^1_0(\Omega)$ into $L^1(\Omega)$. 
The classical energy associated to \eqref{eq} is defined along a trajectory 
$U=(u,\partial_t u)$ as  
$$E(U)= \int_\Omega \frac12(|\grad u|^2+\alpha |u|^2+|\partial_t u|^2) + 
V(x,u)~.$$
The damping effect appears by the computation
\begin{equation}\label{eq-deriv-E}
\partial_t E(U(t))= - \int_\Omega \gamma(x)|\partial_t u |^2~.
\end{equation}
In particular, the energy $E$ is non-increasing along the trajectories. 
Moreover, the sign assumption \eqref{hyp-f-sign} yields that $V(x,u)\geq 0$. 
Thus, we have that $E(U)\geq C \|U\|_X^2$ and that $E(t)$, $t\geq 0$, is 
bounded on the bounded sets of $X$. All together, the above properties show that 
for any $U_0\in X$, the solution $U=(u,\partial_t u)$ of \eqref{eq} is defined 
for all non-negative times and remains in a bounded set of $X$, which only 
depends on $\|U_0\|_X$. 

A fundamental question of this paper concerns the solution for which the energy 
is constant: are they equilibrium points or may they be moving trajectories? 
At least, the answer is known for the linear equation, see \cite{Dafermos} and 
also \cite{Haraux}. 
\begin{theorem}\label{th-stab-linear}
{\bf Dafermos (1978).}\\
Assume that the damping $\gamma\geq 0$ does not vanish everywhere. Then, for 
any $U_0\in X$, we have 
$$e^{At}U_0~\xrightarrow[~~t\longrightarrow +\infty~~]{}~0~\text{ in }X~.$$
\end{theorem}

\section{Asymptotic compactness and reduction to a unique continuation 
problem}\label{section-asymp-compact}

In this section, we assume a fast enough semi-uniform linear decay as 
described by \eqref{hyp-decay-1} and \eqref{hyp-decay-2}.
We first notice that, by linear interpolation, we have the following result.
\begin{prop}\label{prop-decay-interpole}
For any $\sigma_1,\sigma_2$ such that $0\leq \sigma_2<\sigma_1 \leq 1$, the 
linear semigroup is 
well defined from $X^{\sigma_1}$ in $X^{\sigma_2}$ and we have 
$$\forall U_0\in X^{\sigma_1}~,~~ \|e^{At}U_0\|_{X^{\sigma_2}} \leq 
h(t)^{\sigma_1-\sigma_2} \|U_0\|_{X^{\sigma_1}}~.$$
\end{prop}
\begin{demo}
 We interpolate the estimates \eqref{eq-bound-semigroup} for $\sigma=1$ 
and \eqref{hyp-decay-1} with respective weights 
$(\sigma_2/\sigma_1,1-\sigma_2/\sigma_1)$ and obtain
$$\forall U_0\in D(A)~,~\|e^{At}U_0\|_{X^{\sigma_2/\sigma_1}} \leq 
h(t)^{1-\sigma_2/\sigma_1}\|U_0\|_{D(A)}~.$$
It remains to interpolated the above estimate and \eqref{eq-bound-semigroup} 
for $\sigma=0$ with respective weights $(\sigma_1,1-\sigma_1)$. 
\end{demo}

We also need some regularity properties for $F$. The following properties 
depend on Sobolev embeddings and so of the dimension $d$ of $\Omega$. For 
$d=2$, which is the case in our examples, the properties are general. For 
$d=3$, they are more restrictive but they are shown in the same way. We choose 
to also consider this case in our paper for possible later uses.
\begin{prop}\label{prop-compact-F}
Assume that dim$(\Omega)=2$. Then for any $\sigma\in [0,1)$, the function $F$ 
maps any bounded set $\Bc$ of $X=H^1_0(\Omega)\times L^2(\Omega)$ in a bounded 
set $F(\Bc)$ contained in $X^\sigma= (H^{1+\sigma}(\Omega)\cap 
H^1_0(\Omega))\times H^\sigma_0(\Omega)$. Moreover, $F(\Bc)$ has compact 
enclosure in $X^\sigma$.

If dim$(\Omega)=3$ and if \eqref{hyp-f-estim} holds for 
some $p\in [0,3)$, then the same properties hold for $\sigma\in [0,(3-p)/2)$.
\end{prop}
\begin{demo}
Assume that dim$(\Omega)=2$. 
First notice that we only have to show that $f(\cdot,u)$ is compactly bounded 
in $H^\sigma_0(\Omega)$ since the first component of $F$ is zero. Also 
notice that $f(x,0)=0$ due to the sign assumption \eqref{hyp-f-sign}, thus the 
Dirichlet boundary condition possibly contained in $H^\sigma_0(\Omega)$ 
will be fulfilled by $f(x,u)$ if $u\in H^1_0(\Omega)$. Due to the Sobolev 
embeddings, and since $\Omega$ is compact, it is sufficient to show 
that $F(\Bc)$ is bounded in $W^{1,q}(\Omega)$ for all $q\in [1,2)$ to obtain 
compactness in $H^\sigma(\Omega)$ for any $\sigma\in [0,1)$. Since 
$f$ is of polynomial type due to \eqref{hyp-f-estim}, we know that $f(x,u)$ is 
bounded in $L^q(\Omega)$ for any $q\in [1,2)$. On the other 
hand, using \eqref{hyp-f-estim}, we have
\begin{align*}
\|\grad(f(x,u))\|_{L^q} & \leq  \|\grad_xf(x,u)\|_{L^q} + \|f'_u(x,u)\grad 
u\|_{L^q}\\
& \leq C \|(1+|u|)^p \|_{L^q} + C \|(1+|u|)^{p-1} \grad 
u\|_{L^q}\\
& \leq C(1+\|u\|_{L^{pq}}^p + \|\grad u\|_{L^2} \|u\|^{p-1}_{L^r})
\end{align*}
with $r=(p-1)\frac {2q}{2-q}$ defined as soon as $q<2$. This shows that 
$f(\cdot,u)$ belongs to $W^{1,q}(\Omega)$ for any $q\in [1,2)$ and concludes 
the proof for dim$(\Omega)=2$. 

The case dim$(\Omega)=3$ is similar once we use 
the suitable Sobolev embeddings.
\end{demo}

The main results of this section are the following asymptotic compactness 
properties.
\begin{prop}\label{prop-asymp-compact-1}
If dim$(\Omega)=2$, set $\sigma_*=\sigma_h$. If dim$(\Omega)=3$, assume that 
$p\in[0,3)$ in \eqref{hyp-f-estim} and that $\sigma_h\in ((p-1)/2,1)$ in 
\eqref{hyp-decay-2} and set $\sigma_*=\sigma_h - (p-1)/2$.

Let $U(t)=(u,\partial_t u)$ where $u$ solves \eqref{eq} and let $(t_n)$ be a 
sequence of times such that $t_n\rightarrow +\infty$. Then, there exist a 
subsequence $(t_{\varphi(n)})$ and a solution $W(t)=(w,\partial_t w)(t)$ of 
\eqref{eq} defined for all $t\in\RR$, such that
$$\forall 
t\in\RR~,~~U(t_{\varphi(n)}+t)~\xrightarrow[~~n\longrightarrow 
+\infty~~]{}~W(t)~~~\text{ in }X=H^1_0(\Omega)\times L^2(\Omega)~. $$
Moreover, the solution $W$ is globally bounded in $X^\sigma$ for all  
$\sigma\in [0,\sigma_*)$ and the energy $E(W(t))$ is constant.
\end{prop}
\begin{demo}
Assume first that dim$(\Omega)=2$. We have
\begin{equation}\label{eq-prop-asymp-compact-1}
U(t_n)=e^{At_n}U_0 + \int_0^{t_n} e^{A(t_n-s)}F(U(s))\,{\rm d}s~.
\end{equation}
Due to Theorem \ref{th-stab-linear}, the term $e^{At_n}U_0$ goes to zero in 
$X$. Thus, it remains to show that $\int_0^{t_n} e^{A(t_n-s)}F(U(s))\,{\rm d}s$ 
is a compact term in $X$. First notice that $U(s)$ is uniformly bounded for 
$s\geq 0$ due to the non-increasing energy $E$ (see Section 
\ref{section-basic}).
Due to Proposition \ref{prop-compact-F}, $F(U(s))$ thus 
belongs to a bounded set of $X^{\sigma_1}$ for all $\sigma_1\in 
[0,1)$. By Assumption \eqref{hyp-decay-2} and Proposition 
\ref{prop-decay-interpole}, $e^{A(t_n-s)}F(U(s))$ has an integral in 
$[0,t_n]$ bounded in $X^{\sigma_2}$ uniformly with respect to $n$, 
for any $\sigma_2 \in (0,\sigma_h)$. Thus $\int_0^{t_n} 
e^{A(t_n-s)}F(U(s))\,{\rm d}s$ is a compact sequence in $X^\sigma$ for any
$\sigma\in [0,\sigma_2)$. As a consequence, for any $\sigma$ as 
close as wanted to $\sigma_h$, we may extract a subsequence $(t_{\varphi(n)})$ 
such that $\int_0^{t_{\varphi(n)}} e^{A(t_{\varphi(n)}-s)}F(U(s))\,{\rm d}s$
converges to some limit $W(0)$ in $X^\sigma$. Since the linear term of 
\eqref{eq-prop-asymp-compact-1} goes to zero in $X$ for 
$t_{\varphi(n)}\rightarrow +\infty$, $U(t_{\varphi(n)})$ converges to $W(0)\in 
X^\sigma$ for the norm of $X$.

Let $W(t)=(w,\partial_t w)(t)$ be the maximal solution of the damped wave 
equation \eqref{eq} corresponding to the initial data $W(0)$. 
Let $t\in\RR$, for $n$ large enough $t_{\varphi(n)}+t\geq 0$ and 
thus $U(t_{\varphi(n)}+t)$ is well defined and uniformly bounded in $X$.
Since our equation is well posed, the solution is continuous with respect to the 
initial data. Thus, since $U(t_{\varphi(n)})$ converges to $W(0)$ 
in $X$, we have that $U(t_{\varphi(n)}+t)$ converges to $W(t)$ for all 
$t$ such that $W(t)$ is well defined. But due to the uniform bound on 
$U(t_{\varphi(n)}+t)$, $W(t)$ is uniformly bounded and thus the solution may be 
extended to a global solution $W(t)$, $t\in\RR$. In addition, the 
$X^\sigma$-bound obtained above for $W(0)$ only depends on the $X-$bound on 
$U(s)$ which is uniform due to non-increase of 
the energy of $U(t)$. Thus, the same arguments applied to the convergence 
$U(t_{\varphi(n)}+t)\rightarrow W(t)$ give the same $X^\sigma$-bound for $W(t)$ 
for all $t\in\RR$. Finally, since the energy of $U(t)$ is non-increasing and 
non-negative, for any $t\in\RR$, we must have $E(U(t_n+t))-E(U(t_n))\rightarrow 
0$ when $n\rightarrow 0$ (since $t_n$ goes to $+\infty$). This shows that 
$E(W(t))$ is constant and finishes the proof.

The case dim$(\Omega)=3$ is similar once we take into account the constraints 
given by Proposition \ref{prop-compact-F}.
\end{demo}

\begin{prop}\label{prop-asymp-compact-2}
If dim$(\Omega)=2$, set $\sigma_*=\sigma_h$. If dim$(\Omega)=3$, assume that 
$p\in[0,3)$ in \eqref{hyp-f-estim} and that $\sigma_h\in ((p-1)/2,1)$ in 
\eqref{hyp-decay-2} and set $\sigma_*=\sigma_h - (p-1)/2$.

Let $\sigma>0$ and $R>0$. Let $U_n(t)=(u_n,\partial_t u_n)$ a sequence of 
solutions $u_n$ of \eqref{eq} such that $(U_n(0))\subset X^\sigma$ and 
$\|U_n(0)\|_{X^\sigma}\leq R$. Let $(t_n)$ be a 
sequence of times such that $t_n\rightarrow +\infty$ and let $\sigma'\in 
[0,\min(\sigma,\sigma_*))$. Then, there exist subsequences $(t_{\varphi(n)})$ 
and $(U_{\varphi(n)})$ and a solution $W(t)=(w,\partial_t w)(t)$ of 
\eqref{eq} defined for all $t\in\RR$, such that
$$\forall 
t\in\RR~,~~U_{\varphi(n)}(t_{\varphi(n)}+t)~\xrightarrow[~~n\longrightarrow 
+\infty~~]{}~W(t)~~~\text{ in } X^{\sigma'}~. $$
Moreover, the solution $W$ is globally bounded in $X^{\sigma'}$ and the energy 
$E(W(t))$ is constant.
\end{prop}
\begin{demo}
The arguments are similar as the ones of the above proof of Proposition 
\ref{prop-asymp-compact-1}. The term $e^{At_n}U_n(0)$ goes to zero in 
$X^{\sigma'}$ due to Proposition \ref{prop-decay-interpole} because 
$\sigma'<\sigma$. We bound the integral $\int_0^{t_n} 
e^{A(t_n-s)}F(U_n(s))\,{\rm d}s$ as in the proof of Proposition 
\ref{prop-asymp-compact-1}: $U_n(s)$ is uniformly bounded in $X$, so 
$F(U_n(s))$ is uniformly bounded in $X^\eta$ with $\eta<1$ in dimension $2$ or 
$\eta<(3-p)/2$ in dimension $3$. Proposition \ref{prop-decay-interpole} 
together with \eqref{hyp-decay-2} implies that the integral is uniformly bounded 
in $X^{\sigma'}$ with $\sigma'<\sigma_*$. The compactness follows by leaving 
any small amount of regularity in the process. To obtain the convergence to 
$W(t)$ for all $t$, we use the same argument as the one of the proof of 
Proposition \ref{prop-asymp-compact-1} to first show the convergence in $X$. 
Then, the above arguments also show the compactness of $U(t_{\varphi(n)}+t)$ in 
$X^{\sigma'}$ and thus the convergence to $W(t)$ also holds in 
$X^{\sigma'}$. The last property is the same as the ones 
of Proposition \ref{prop-asymp-compact-1}. 
\end{demo}

The conclusions of Theorem \ref{th1} then follow from Propositions 
\ref{prop-asymp-compact-1} and 
\ref{prop-asymp-compact-2} as soon as we can prove that $W\equiv 0$ for any 
subsequences of any sequences $(t_n)$ and $U_n$. To this end, notice that 
$E(W(t))$ is constant and its derivative \eqref{eq-deriv-E} implies that 
$\int_\Omega \gamma(x) |\partial_t w|^2=0$ for all time. To formulate this 
property as a unique continuation property, we set as usual $z=\partial_t w$ 
and notice that $z$ solves
\begin{equation}\label{eq-UCP}
\left\{\begin{array}{ll}
\partial_{tt}^2 z(x,t) = \Delta z(x,t) - \alpha z(x,t) - 
f'_u(x,w(x,t)) z~~~&(x,t)\in\Omega\times\RR\\
z_{|\partial\Omega}(x,t)=0&(x,t)\in\partial\Omega\times \RR\\
z(x,t)\equiv 0 & (x,t)\in \text{support}(\gamma)\times \RR
\end{array}\right.
\end{equation}
If this implies $z\equiv 0$ everywhere, this means that $w(x,t)=w(x)$ is 
constant in time and solves
$$\Delta w(x)-\alpha w(x)=f(x,w(x))~.$$
Multiplying by $w$ and integrating, we obtain 
$$\|\grad w\|^2 +\alpha \|w\|^2 = -\int_\Omega f(x,w(x))w(x)\,{\rm d}x~.$$
By the sign Assumption \eqref{hyp-f-sign}, this yields $w\equiv 0$. Thus, it 
only remains to study this unique continuation property.
\begin{prop}\label{prop-UCP-implique-th}
Assume that $z\equiv 0$ is the only global solution of \eqref{eq-UCP}. Then the 
decay assumptions \eqref{hyp-decay-1} and \eqref{hyp-decay-2} imply the 
conclusions of Theorem \ref{th1}.
\end{prop}

\section{Rate of the nonlinear decay: proof of Proposition 
\ref{prop1}}\label{section-decay}
The purpose of this section is to prove Proposition \ref{prop1}.
When estimating the decay rate of the nonlinear system, we will not exactly 
need the decay of the linear semigroup but more precisely the decay rate of the 
linearization at $u=0$. This is not difficult since an estimate as 
\eqref{hyp-decay-1} is a high-frequency result: the behavior of the high 
frequencies is the difficult part and we only need that the low frequencies do 
not lie on the imaginary axes.

\subsection{The polynomial case}
The case of polynomial decay is obtained as follows.
\begin{lemma}\label{lemma-appli1}
Assume the sign hypothesis \eqref{hyp-f-sign} and assume that 
\eqref{hyp-decay-1} holds with $h(t)=\Oc(t^{-\alpha})$, with $\alpha>1$. Set 
$$\tilde A=A + \left(\begin{array}{c} 0 \\ 
-f'_u(x,0)\end{array}\right)=\left(\begin{array}{cc} 0 & Id\\ 
\Delta-\alpha Id-f'_u(x,0) & -\gamma(x) \end{array}\right)~.$$
Then, there exists $C>0$ such that 
$$\forall t\geq 0~,~\forall U_0\in D(A)~,~~\|e^{\tilde At} U_0\|_X \leq \frac 
C{(1+t)^{\alpha}} \|U_0\|_{D(A)}~.$$
\end{lemma}
\begin{demo}
Due to \eqref{hyp-f-sign}, we have that $f(x,0)=0$ and that $f'_u(x,0)\geq 0$. 
Since $\tilde A$ is a compact perturbation of $A$, we do not expect that the 
behavior for high frequencies should be modified. For low frequencies, the sign 
of $f'_u(x,0)$ is sufficient to avoid eigenvalues on the imaginary axes. 

To prove rigorously these facts as quickly as possible, we use the results 
stated in Appendix with $H=L^2(\Omega)$, 
$L=-\Delta+\alpha$, $B=\gamma$ and $V=f'_u(x,0)$. Due to Theorem \ref{th-BT}, 
there exist $\mu_0$ and $C>0$ such that 
$$\forall \mu\in\RR\text{ with }|\mu|\geq 
\mu_0~,~~\|(A-i\mu)^{-1}\|_{\Lc(L^2)} 
\leq C {|\mu|^{1/\alpha}}~.$$
We now use Proposition \ref{prop-appendix-linearized} to obtain that the same 
estimate holds for the resolvents $(\tilde A-i\mu)^{-1}$ for large 
$\mu$. Moreover, for $\mu$ in a compact interval, Proposition 
\ref{prop-resolvent} ensures that the resolvent $(\tilde A-i\mu)^{-1}$ is well 
defined. Applying Theorem \ref{th-BT} in the converse sense concludes the proof.
\end{demo}

{ \noindent \emph{\textbf{Proof of the first case of Proposition \ref{prop1}:}}}
we assume that the conclusions of Theorem \ref{th1} hold. In particular, 
the trajectory of a ball of $X^\sigma$ of radius $R$ is attracted by $\{0\}$ in 
$X^{\sigma'}$ for a small enough $\sigma'>0$. Thus, it is sufficient to 
prove that the decay has the same rate as the linear one, as soon as we start 
from a small ball of $X^{\sigma'}$ of radius $\rho>0$ and stay in it. 

We consider the linearization of our equation near the stable state $u=0$. We 
set $\tilde A$ be as in Lemma \ref{lemma-appli1} and $\tilde 
F(u)=(0,f(x,u)-f'_x(x,u)u)$. 
Since $H^{1+\sigma'}(\Omega)$ is embedded in $L^\infty(\Omega)$ and is an 
algebra, we may bound the derivatives of $f$ and by linearization, for any 
small $\delta>0$, we may work with $U(t)$ in a ball of $X^{\sigma'}$ of radius 
$\rho$, which is such that $\|\tilde F(U)\|_{X^1} \leq \delta \|U\|_{X}$.

Let $U(t)$ be a trajectory in the small ball of 
$X^{\sigma'}$ with $\|U_0\|_{X^\sigma}\leq \rho$. We have
$$(1+t)^{\sigma\alpha} U(t)=(1+t)^{\sigma\alpha} e^{\tilde At}U_0 + 
(1+t)^{\sigma\alpha} \int_0^t e^{\tilde A(t-s)}\tilde F(U(s))\,{\rm d}s~.$$
The term $(1+t)^{\sigma\alpha} e^{\tilde At}U_0$ is bounded by the linear decay 
(see Lemma \ref{lemma-appli1} and Proposition \ref{prop-decay-interpole}). By 
using the above estimate on $\tilde F(U(s))$, we get that
$$\|(1+t)^{\sigma\alpha} U(t)\|_X \leq C + (1+t)^{\sigma\alpha} \int_0^t 
\frac{C}{(1+(t-s))^\alpha}\delta \|U(s)\|_{X}\,{\rm d}s~.$$
Thus, 
\begin{align*}
\max_{t\in[0,T]} (1+t)^{\sigma\alpha} \|U(t)\|_X ~\leq~& C + \delta C 
\left(\max_{s\in [0,T]} 
(1+s)^{\sigma\alpha} \|U(s)\|_X \right) \\ & ~~~\times \max_{t\in[0,T]} 
\int_0^t 
\left(\frac{1+t}{1+s}\right)^{\sigma\alpha} \frac{{\rm d}s}{(1+(t-s))^\alpha}
\end{align*}
where $C$ is a constant independent on $T$ when $T$ goes to $+\infty$ and on 
the radius $\rho$ of the starting ball when $\rho$ goes to $0$. The limit 
$T\rightarrow +\infty$ 
will prove our theorem as soon as we can show that the integral term is bounded 
uniformly in $t$. Indeed, up to work with $\rho$ small 
enough, we may assume that $\delta$ is such that the whole last term is less 
than $1/2 \max_{s\in [0,T]} (1+s)^{\sigma\alpha} \|U(s)\|_X$ and may be 
absorbed by the left hand side.

To estimate the integral, we use the change of variable $\tau = (1+s)/(t+2)$, 
for which $1+(t-s)=(t+2)(1-\tau)$. We obtain that 
$$I(t)=\int_0^t \left(\frac{1+t}{1+s}\right)^{\sigma\alpha} 
\frac{{\rm d}s}{(1+(t-s))^\alpha} = \left(\frac{1+t}{2+t}\right)^{\sigma\alpha} 
\frac{1}{(t+2)^{\alpha-1}} \int_{1/(t+2)}^{1-1/(t+2)} \frac {{\rm 
d}\tau}{\tau^{\sigma\alpha}(1-\tau)^\alpha}~.$$
Recall that $\alpha>1$ and that $\sigma \leq 1$. The integral 
$\int_0^1 \frac {{\rm d}\tau}{\tau^{\sigma\alpha}(1-\tau)^\alpha}$
does not converge at least close to $1$ and if it also diverges close to $0$, 
the blow up is slower or equal to the one occurring close to $1$. Thus, 
$\int_{1/(t+2)}^{1-1/(t+2)} \frac {{\rm 
d}\tau}{\tau^{\sigma\alpha}(1-\tau)^\alpha}$ is of order $\Oc(t^{\alpha-1})$ 
when $t$ goes to $+\infty$. This shows that the whole integral $I(t)$ is bounded 
uniformly in $t$. 
{\hfill$\square$\\}

\subsection{The exponential case}
The following lemma is similar to Lemma \ref{lemma-appli1}, except that we 
cannot use the result of Borichev and Tomilov recalled in Theorem \ref{th-BT}. 
If the decay is not polynomial, then we must accept a logarithmic loss and use 
the results of Batty and Duyckaerts, Theorems \ref{th-BD-1} and \ref{th-BD}.
\begin{lemma}\label{lemma-appli2}
Assume the sign hypothesis \eqref{hyp-f-sign} and assume that 
\eqref{hyp-decay-1} holds with $h(t)=\Oc(e^{-a t^{1/\beta}})$, with 
$a>0$ and $\beta>0$. Set 
$$\tilde A=A + \left(\begin{array}{c} 0 \\ 
-f'_u(x,0)\end{array}\right)=\left(\begin{array}{cc} 0 & Id\\ 
\Delta-\alpha Id-f'_u(x,0) & -\gamma(x) \end{array}\right)~.$$
Then, there exists $C>0$ and $b>0$ such that 
$$\forall t\geq 0~,~\forall U_0\in D(A)~,~~\|e^{\tilde At} U_0\|_X \leq 
Ce^{-b t^{1/(\beta+1)}} \|U_0\|_{D(A)}~.$$
\end{lemma}
\begin{demo}
As in the proof of Lemma \ref{lemma-appli1}, we use the results 
stated in Appendix with $H=L^2(\Omega)$, 
$L=-\Delta+\alpha$, $B=\gamma$ and $V=f'_u(x,0)$. 
Using Theorem \ref{th-BD-1}, we obtain that
$$\forall \mu\in\RR\text{ with }|\mu|\geq 
\mu_0~,~~\|(A-i\mu)^{-1}\|_{\Lc(L^2)} \leq C {|\ln \mu|^{\beta}}~.$$
As in Lemma \ref{lemma-appli1}, we use Proposition 
\ref{prop-appendix-linearized} to obtain that the same 
estimate holds for the resolvents $(\tilde A-i\mu)^{-1}$ for large 
$\mu$ and Proposition \ref{prop-resolvent} to deal with the low frequencies.
The difference is that Theorem \ref{th-BD} yields a logarithmic loss when going 
back to the estimate of the semigroup (see the definition of $M_{\text{log}}$), 
leading to the exponent $t^{1/(\beta+1)}$.
\end{demo}

\noindent {\bf Remark:} We have seen that there is a logarithmic loss in our 
estimate. However, in the applications, we will obtain a better result. Indeed, 
this loss was already present in the original estimate for the linear semigroup 
because of the additional log in $M_{\text{log}}$ of Theorem \ref{th-BD}. In 
some sense, the above abstract result makes an additional use of the back and 
forth Theorems \ref{th-BD-1} and \ref{th-BD}. We can improve our estimate by a 
shortcut: we go back to the estimate of the resolvent in the original proof 
of the linear decay, before the authors apply Theorem \ref{th-BD}, and 
we directly apply the above arguments to estimate $(\tilde A-i\mu)^{-1}$ and 
then apply Theorem \ref{th-BD}. With this trick, we do not add a second 
logarithmic loss to the one of the original proof dealing with the linear 
semigroup. However, we can do this only in the concrete situations and 
not in an abstract result as Proposition \ref{prop1}.\\[3mm]

{ \noindent \emph{\textbf{Proof of the  second case of Proposition 
\ref{prop1}:}}} the method is exactly the same as in the first case. The only 
difference is that, instead of bounding $\int_0^t 
\left(\frac{1+t}{1+s}\right)^{\sigma\alpha} 
\frac{{\rm d}s}{(1+(t-s))^\alpha}$, we must here bound an integral of the type
$$I(t)=\int_0^t e^{\sigma c(t^\gamma-s^\gamma)} 
e^{-c (t-s)^\gamma}\,{\rm d}s$$
for some $c>0$ and $\gamma \in(0,1)$. We set $\tau=s/t$ and obtain
$$I(t)=t\int_0^1 
e^{c{t^\gamma}(\sigma-\sigma \tau^\gamma-(1-\tau)^\gamma)}\,{\rm 
d}\tau \leq t\int_0^1 
e^{c \sigma t^\gamma(1-\tau^\gamma-(1-\tau)^\gamma)}\,{\rm 
d}\tau$$
and by symmetry
$$I(t)\leq 2t\int_0^{1/2}
e^{ \sigma c t^\gamma(1-{\tau}^\gamma-(1-\tau)^\gamma)}\,{\rm 
d}\tau~.$$
We notice that $\tau\mapsto 1-\tau^\gamma-(1-\tau)^\gamma$ is decreasing for 
$\tau\in [0,1/2]$ since its derivative is 
$\gamma((1-\tau)^{\gamma-1}-\tau^{\gamma-1})$ with $\gamma-1 < 0$. Moreover, 
$1-\tau^\gamma-(1-\tau)^\gamma \sim -\tau^\gamma$ for small $\tau$. Thus, there 
exists $\nu>0$ small enough such that $1-\tau^\gamma-(1-\tau)^\gamma \leq 
-\nu \tau^\gamma$ for $\tau\in[0,1/2]$. We get
$$
I(t)\leq 2t\int_0^{1/2} e^{\sigma c t^\gamma (-\nu 
\tau^\gamma)}\,{\rm d}\tau= 2\int_0^{t/2} e^{-\sigma c \nu s^\gamma}\,{\rm d}s~.
$$
The integrand of these last bound is integrable on $\RR_+$, thus $I(t)$ is 
bounded uniformly with respect to $t$. Arguing as in the proof of the 
polynomial case, this proves the second part of Proposition \ref{prop1}.
{\hfill$\square$\\}

\section{Application 1: the open book}\label{section-appli1}

In this section, we consider the third example of Figure \ref{fig-1}. Let 
$\Omega=\TT^2$ be the two-dimensional torus and let $\alpha>0$ (there is no 
boundary and so no Dirichlet boundary condition). Assume that 
$$\gamma(x_1,x_2)=|x_1|^\beta~.$$
We have the following decay estimate proved in \cite[Theorem 
1.7]{Leautaud-Lerner}.
\begin{theorem}\label{th-Matthieu-Lerner}
{\bf L\'eautaud \& Lerner, 2015.}\\
In the above setting, the semigroup $e^{At}$ satisfies
$$\|e^{At} U_0\|_X \leq \frac C{(1+t)^{1+2/\beta}} \|U_0\|_{D(A)}~.$$
\end{theorem}
Of course, this estimate implies the decay assumptions 
\eqref{hyp-decay-1} and \eqref{hyp-decay-2} for $\sigma_h < 2/(2+\beta)$. Since 
the support of $\gamma$ is $\TT^2$, the unique continuation property is trivial 
and Proposition \ref{prop-UCP-implique-th} implies that the conclusions 
of Theorem \ref{th1} holds in this case. Moreover, Proposition \ref{prop1} 
provides an explicit decay rate, which is optimal (since it is the same as the 
linear one). Due to the trivial unique continuation property, this case is 
far simpler than the general results Theorem \ref{th1} and \ref{th2}. 
Nevertheless, it seems the first non-linear stabilization and decay 
estimate in a case where the linear semigroup has only a polynomial decay. 
\begin{theorem}\label{th-appli1}
Consider the damped wave equation \eqref{eq} in $\Omega=\TT^2$ and with 
$\gamma(x_1,x_2)=|x_1|^\beta$ ($\beta>0$), $\alpha>0$ and $f$ satisfying 
\eqref{hyp-f-estim} and \eqref{hyp-f-sign}. Then, any solution $u$ of 
\eqref{eq} satisfies 
$$\|(u,\partial_t u)(t)\|_{H^1_0\times L^2}  ~\xrightarrow[~~t\longrightarrow 
+\infty~~]{}~0~.$$

Moreover, for any $R$ and $\sigma\in (0,1]$, there exists 
$C_{R,\sigma}$ such that, for any solution $u$ with $U_0\in 
X^\sigma$, 
$$ \|(u_0,u_1)\|_{H^{1+\sigma}\times H^\sigma}\leq 
R~~\Longrightarrow ~~ \|(u,\partial_t u)(t)\|_{H^1_0\times L^2} \leq 
\frac{C_{R,\sigma}}{(1+t)^{\sigma(1+2/\beta)}}~.$$
\end{theorem}


\section{Unique continuation theorems}\label{section-UCP}
As proved in Section \ref{section-asymp-compact}, the last step to prove 
stabilization is the unique continuation property: if $z$ is a global solution 
of  
\begin{equation}\label{eq-UCP2}
\left\{\begin{array}{ll}
\partial_{tt}^2 z(x,t) = \Delta z(x,t) - \alpha z(x,t) - 
f'_u(x,w(x,t)) z~~~&(x,t)\in\Omega\times\RR\\
z_{|\partial\Omega}(x,t)=0&(x,t)\in\partial\Omega\times \RR\\
z(x,t)\equiv 0 & (x,t)\in \text{support}(\gamma)\times \RR
\end{array}\right.
\end{equation}
then $z\equiv 0$. Except for the example of Section \ref{section-appli1}, the 
property is often difficult to obtain. The purpose of this section is to gather 
several results yielding this property.

The first known result has been proved by Ruiz in \cite{Ruiz}. It stated the 
unique continuation property in a bounded domain $\Omega \subset\RR^d$ as 
soon as the support of $\gamma$ contains a neighborhood of the boundary 
$\partial\Omega$. This result has been generalized in \cite{K-K} (see 
also \cite{Las-Tri-Zha} for Neumann boundary conditions).
However, this kind of results is not relevant in this paper. Indeed, their 
geometric settings implies the uniform decay of the semigroup $e^{At}$ and we 
are interested here in cases where it is not satisfied.
We need sharper results.

\subsection{Unique continuation with coefficients analytic in time}

A very general unique continuation property holds if the 
coefficients of a linear wave equation as \eqref{eq-UCP2} are analytic in time.
This is a consequence of local continuation results proved by 
by H\"ormander in \cite{Hormander_uniq_analytic} and generalized by Tataru in 
\cite{Tataru_uniq_analyticJMPA} and also independently proved by 
Robbiano and Zuily in \cite{Rob-Zui}. These results concern 
in fact a very general setting but we restrict here the statements at the case 
of the wave equation. The application to the wave equation and the proof that 
the local results yield a global one are classical and straightforward, see for 
example \cite[Corollary 3.2]{RC} for the details.
\begin{theorem}\label{th-RZ}
{\bf Robbiano-Zuily, H\"ormander (1998)}\\
Let $T>0$ (or $T=+\infty$) and let $b$, $c$ and
$d$ be smooth coefficients. Assume moreover that $b$, $c$ and $d$ are analytic 
in time and that $z$ is a strong solution of
\begin{equation}\label{eq-RZ}
\left\{\begin{array}{ll}
\partial^2_{tt} z=\Delta z + b(x,t)\partial_t z+c(x,t).\grad
z+d(x,t)z~~~~~&(x,t)\in\Omega\times (-T,T)\\
z_{|\partial\Omega}(x,t)=0&(x,t)\in\partial\Omega\times \RR~.
\end{array}\right.
\end{equation}
Let $\Oc$ be a non-empty open subset of $\Omega$ and assume that 
$z(x,t)=0$ in $\Oc\times (-T,T)$. Then
$z(x,0)\equiv 0$ in $\Oc_T=\{x_0\in\Omega~,~d(x_0,\Oc)<T\}$.

As consequences if $z\equiv 0$ in $\Oc\times (-T,T)$ and $\overline 
\Oc_T=\Omega$, then $z\equiv 0$ everywhere.
\end{theorem}

\subsection{Unique continuation through pseudo-convex surfaces without boundary}
\label{subsectUCPwithout}
If the coefficients of \eqref{eq-RZ} are not analytic in time, the geometry of 
the problem is more constrained. However, it could still include cases where 
the geometric control condition of \cite{BLR} does not hold and thus where the 
semigroup $e^{At}$ is not uniformly stable, see the examples below.

We consider here H\"ormander framework (see \cite{Hormander} for example). The 
principal symbol of the differential operator of \eqref{eq-UCP2} is of order 
two and writes locally
$$p(x,t,\xi,\tau)=\xi\trans.A(x).\xi-|\tau|^2$$
where $A(x)$ is a smooth family of positive definite symmetric matrices coding 
the Beltrami Laplacian operator in a local chart. Let $\phi(x,t,\xi,\tau)$ be 
a locally $\Cc^1-$function, we introduce the Poisson bracket 
$$H_p(\psi)=\{p,\psi\}= \grad_\xi p \grad_x \psi + \partial_\tau p \partial_t 
\psi - \grad_x p \grad_\xi \psi - \partial_t p \partial_\tau \psi~.$$
Let $\psi$ be a smooth function defined in a neighborhood $\Oc\subset 
\RR^{d+1}$ of $(x_0,t_0)$. Assume that $(\grad_x \psi,\partial_t 
\psi)(x_0,t_0)\neq 0$ so that $\Sigma=\{(x,t),\psi(x,t)=0\}$ defines a 
smooth hypersurface near $(x_0,t_0)$.
\begin{defi}\label{defi-pseudoconvex}
The local hypersurface $\Sigma$ is said to be \emph{non-characteristic} at 
$(x_0,t_0)$ if 
$$p(x_0,t_0,\grad \psi(x_0,t_0),\partial_t \psi(x_0,t_0))\neq 0~.$$
Moreover, $\Sigma$ is said to be \emph{strongly pseudo-convex} at $(x_0,t_0)$ 
if for any $(\xi,\tau)\neq 0$ such that $p(x_0,t_0,\xi,\tau)=0$ and 
$H_p(\psi)(x_0,t_0)=0$, we have 
$$H_p^2(\psi)(x_0,t_0)>0~.$$   
\end{defi}

Notice that the above definition of strongly pseudo-convexity is adapted to the 
case of a real differential operator of order two. Thus it is perfectly adapted 
to the situation of this paper where the wave operator is 
$p=\xi\trans.A(x).\xi-|\tau|^2$. However, we emphasize that, in the general 
case, the assumption of pseudo-convexity is more complex, see \cite{Hormander}.

The geometrical interpretation of Definition \ref{defi-pseudoconvex} is as 
follows. First, the fact that the surface is non-characteristic says that 
$|\partial_t \psi|^2 \neq \grad\psi\trans.A(x).\grad \psi$. This means that the 
surface is not moving at the exact same speed as the sound waves. 

The pseudo-convexity is slightly more involved. Consider the total Hamiltonian 
flow $\sigma\mapsto\varphi_\sigma$ defined by
$$\varphi_0(x,t,\xi,\tau)=(x,t,\xi,\tau)~~~\partial_\sigma 
\varphi_\sigma(x,t,\xi,\tau)=(\grad_\xi p,\partial_\tau 
p,-\grad_x p,-\partial_t p)(\varphi_\sigma)~.$$
Since $p$ is independent of $t$, $\tau$ is constant and thus 
$t(\sigma)=t-2\tau\sigma$, meaning that $\sigma$ is a simple new 
parametrization 
of time. Moreover, $\grad_\xi p$ and $\grad_x p$ are independent of $t$ and 
$\tau$. Thus, $(x,\xi)(\sigma)$ follows the geodesic flow
$$\partial_\sigma (x,\xi)=(\grad_\xi g,-\grad_x g)(x,\xi)$$
where $g(x,\xi)=\xi\trans.A(x).\xi$ is the symbol of the local metric. Assume 
that $p(x,t,\xi,\tau)=0$ at $\sigma=0$. The Hamiltonian being conserved, we 
always have $p(x,t,\xi,\tau)=0$ and $|\tau|^2=\xi\trans A(x)\xi$ is constant in 
$\sigma$: the point $x(\sigma)$ is moving along a geodesic of the metric at 
a speed which is of constant norm $|\tau|$ with respect to the metric.
Let $h$ be a function of $(x,t,\xi,\tau)$, then the Poisson bracket $\{p,h\}$ 
is the derivative at $\sigma=0$ of $h(\varphi_\sigma(x,t,\xi,\tau))$. Thus 
$$\{p,h\}=\{g,h\} -2\tau \partial_t h=\{g,h\} + \partial_\sigma h$$
where $\{g,h\}$ is the derivative along the geodesic $(x,\xi)(\sigma)$ of the 
metric starting at $x$ with speed $\xi$. Thus, the strongly 
pseudo-convexity condition $H_p\psi=0\Rightarrow H_p^2 \psi>0$ means that if a 
geodesic of the surface is tangent to $\Sigma$ in the space-time sense, then it 
must be contained in a non-degenerated sense in the half-space $\psi(x,t)>0$ 
for $t\neq0$. 
Finally notice that if $\psi$ does not depend on time $t$ (as in Definitions 
\ref{defi-foliation1} and \ref{defi-foliation2}), then the strongly 
pseudo-convexity is a classical strong convexity: if a classical geodesic of 
the metric $g$ is tangent to the surface $\psi(x)=0$, it must be 
contained in the half-space $\psi(x)>0$.

Theorem 28.4.3 of \cite{Hormander} mainly comes from \cite{LR} and is stated as 
follows.
\begin{theorem}\label{th-hormander}
{\bf Lerner and Robbiano (1985), H\"ormander.}\\
Let $\Oc$ be a small open neighborhood of a point $(x_0,t_0)$ in $\RR^d\times 
\RR$ and let $A(x)$ be a smooth family of positive definite symmetric matrices. 
Let $b$, $c$ and $d$ be bounded  coefficients. Assume that $z$ is a mild 
solution of
\begin{equation}\label{eq-LR}
\partial^2_{tt} z=\div A(x)\grad z + b(x,t)\partial_t z+c(x,t).\grad
z+d(x,t)z~~~~~(x,t)\in\Oc~.
\end{equation}
Let $\Sigma=\{(x,t),\psi(x,t)=0\}$ be a smooth surface containing $(x_0,t_0)$ 
which is non-cha\-rac\-te\-ris\-tic and strongly pseudo-convex in the sense of 
Definition \ref{defi-pseudoconvex}.

Then, if $u(x,t)=0$ for all $(x,t)\in\Oc$ such that $\psi(x,t)\geq 0$, we have 
$u(x,t)\equiv 0$ in a neighborhood of $(x_0,t_0)$.
\end{theorem}

\begin{figure}[ht]
\begin{center}
\resizebox{6cm}{!}{\input{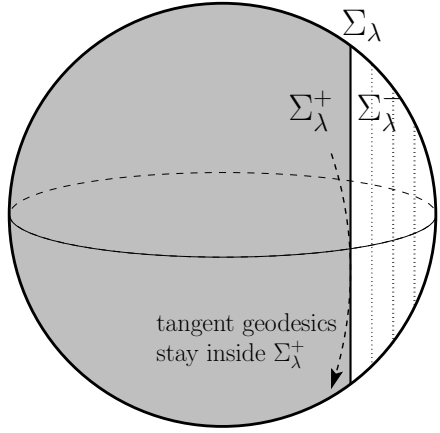}}
\end{center}
\caption{\it A oriented pseudo-convex foliation without boundary in the sphere 
$\SS^2$. The surfaces $\Sigma_\lambda$ forms a smooth family of vertical 
circles inside a hemisphere and $\Sigma_\lambda^+$ and $\Sigma_\lambda^-$ are 
respectively the large and the small spherical caps. The geodesics beings the 
great circles, the ones which are tangent to a surface $\Sigma_\lambda$ stay 
inside $\Sigma_\lambda^+$.\label{fig-foliation-1}}
\end{figure}

Theorem \ref{th-hormander} states a local unique continuation property through 
pseudo-convex surfaces. To use it, it is more convenient to have a global 
version. This kind of global foliation has 
already been introduced in \cite{Stefanov-Uhlmann} by Stefanov and Uhlmann.
\begin{defi}\label{defi-foliation1}
A family of surfaces $(\Sigma_\lambda)_{\lambda\in[0,1)}$ is an \emph{oriented 
pseudo-convex foliation without boundary} in a compact manifold $\Omega$ if:
\begin{enumerate}[(i)]
\item the family of surfaces is smooth in the sense that it is locally described 
as level sets $\{x,\psi_\lambda(x)=0\}$ where 
$(x,\lambda)\mapsto\psi_\lambda(x)$ is a local smooth function with 
$\grad_x\psi_\lambda\neq 0$.
\item each surface is globally oriented in the sense that 
there exist disjoint sets $\Sigma_\lambda^\pm$ such 
that locally $\{x\in\Omega,\pm\psi_\lambda(x)>
0\}\subset \Sigma_\lambda^\pm$ and such that $\Omega=\Sigma_\lambda^- \cup 
\Sigma_\lambda \cup \Sigma_\lambda^+$.
\item for each $\lambda$, $(x,t)\mapsto \psi_\lambda(x)$ is pseudo-convex in the 
sense of Definition \ref{defi-pseudoconvex} as a function independent of $t$. 
Equivalently, $\Sigma_\lambda^-$ is locally strictly convex in a neighborhood of 
its boundary $\Sigma_\lambda$ for the metric $g$: 
for each $x\in\Sigma_\lambda$, a geodesic through $x$ which is tangent at 
$\Sigma_\lambda$ is locally included in $\Sigma_\lambda^+$, $x$ excepted.
\item the surfaces $\Sigma_\lambda$ are compact and have no boundary or 
equivalently do not meet $\partial\Omega$.
\end{enumerate}
\end{defi}
A typical example of such oriented pseudo-convex foliation without boundary is 
given in Figure \ref{fig-foliation-1}.

By a classical argument, we may state a global version of Theorem 
\ref{th-hormander} as follows.  
\begin{theorem}\label{th-hormander-global}
Let $\Omega$ be a smooth compact manifold (with or without boundaries) and let 
$\omega\subset\Omega$ be an open set. Assume that there exists an oriented 
pseudo-convex foliation without boundary $(\Sigma_\lambda)_{\lambda\in[0,1)}$ 
of $\Omega$ in the sense of Definition \ref{defi-foliation1}. Also assume that 
$\Sigma_0^+\subset \omega$ and $\omega\cup \left(\bigcup_{\lambda\in[0,1)} 
\Sigma_\lambda^+\right)$ covers $\Omega$ up to a set of zero measure.

Let $b$, $c$ and $d$ be bounded  coefficients. Assume that $z$ is a global mild 
solution of
$$
\left\{\begin{array}{ll}
\partial^2_{tt} z=\Delta z + b(x,t)\partial_t z+c(x,t).\grad
z+d(x,t)z~~~~~&(x,t)\in\Omega\times \RR\\
z_{|\omega}(x,t)=0&(x,t)\in\omega\times\RR
\end{array}
\right.
$$
with any suitable boundary conditions on $\partial\Omega$ such that the wave 
equation is well-posed and where $\Delta$ is the Laplace-Beltrami operator 
related to $\Omega$.

Then $z\equiv 0$ everywhere.
\end{theorem}
\begin{demo}
We will show that $z(\cdot,t=0)$ vanishes in $\omega\cup\Sigma^+_\lambda$, 
for all $\lambda\in [0,1)$, which shows that $z(t=0)\equiv 0$ and thus that 
$z\equiv 0$ due to the uniqueness  properties of the wave equation. By 
assumption, $z\equiv 0$ in $\Sigma_0^+\subset \omega$. Let $\lambda_0\in (0,1)$ 
and let 
$h_{\alpha,T}(t)=\alpha(1-t^2/T^2)$. We consider the family of surfaces $t\in 
[-T,T]\mapsto \Sigma_{h_{\alpha,T}(t)}$ which is locally parametrized by 
functions $(x,t)\mapsto \psi_{h_{\alpha,T}(t)}(x)$. 
Notice that it is a smooth family of smooth surfaces since due to 
Assumption (i) of Definition \ref{defi-foliation1}. Also notice that the 
larger is $T$, the 
smaller are the derivatives of these functions with respect to $t$. By 
assumption, each function $x\mapsto \psi_\lambda(x)$ is non-characteristic and 
strongly pseudo-convex as a function independent of $t$. By compactness, 
there exists $T$ large enough such that $(x,t)\mapsto \psi_{h_{\alpha,T}(t)}(x)$ 
defines local surfaces which are non-characteristic and pseudo-convex 
for all $\alpha\in[0,\lambda_0]$ and $t\in[-T,T]$. The parameter $T$ is fixed 
in the remaining part of the proof and we may omit it in the notations.

Notice that, for any $\alpha$, the family of set $t\in[-T,T]\mapsto 
\Sigma^+_{h_\alpha(t)}$ starts inside $\omega$ at $t=-T$ and finishes inside 
$\omega$ at $t=T$. Moreover, for any small 
$\alpha$, these sets always stay inside $\omega$ where $z$ vanishes. Assume 
that 
there exist $(x,t)$ and $\alpha\in[0,\lambda_0]$ such that $x\in 
\Sigma^+_{h_\alpha(t)}$ and $z(x,t)\neq 0$. We set 
\begin{equation}\label{eq-demo-global}
\alpha_0=\min\{\alpha\in (0,\lambda_0]~,~
\exists t\in[-T,T],\exists x\in\Sigma^+_{h_{\alpha}(t)}\text{ such that 
}z(x,t)\neq 0\}~.
\end{equation}
By continuity, we know that $z(x,t)=0$ for all $x\in 
\Sigma^+_{h_{\alpha_0}(t)}$. Moreover, there exists $t_0\in (-T,T)$ and 
$x_0\in\Sigma_{h_{\alpha_0}(t_0)}$ such that $z$ is not identically zero in any 
neighborhood of $(x_0,t_0)$. Indeed, otherwise, by compactness, we may extend 
the set where $z$ vanishes and contradict \eqref{eq-demo-global}.

To conclude, it remains to use the local unique continuation property of 
Theorem 
\ref{th-hormander} at $(x_0,t_0)$ with the time-space surface defined by 
$(x,t)\mapsto \psi_{h_{\alpha_0}(t)}(x)$. The continuation implies that $z$ 
vanishes near $(x_0,t_0)$ which contradicts the construction. Thus, $z(x,t)=0$ 
for all $t\in [-T,T]$ and $x\in \bigcup_\alpha \Sigma^+_{h_\alpha(t)}$. In 
particular $z(\cdot,t=0)\equiv 0$ in 
$\bigcup_{\alpha}\Sigma^+_{h_\alpha(0)}=\bigcup_{\lambda\leq\lambda_0}
\Sigma^+_\lambda$. Since these arguments hold for all $\lambda_0<1$ and since 
$\omega\cup_{\lambda\in [0,1)} \Sigma^+_\lambda$ is $\Omega$ up to a set of 
measure zero, we have that $z(\cdot,t=0)\equiv 0$ in $\Omega$. Well-posedness 
of 
the linear wave equation concludes that $z\equiv 0$ everywhere.
\end{demo}

Notice that, as it is stated, this unique continuation result needs an infinite 
time to be efficient, where Theorem \ref{th-RZ} only need a finite explicit 
time. In fact, a careful look to the proof shows that a finite time is 
sufficient once we know that the family of surface is pseudo-convex and 
non-characteristic in a uniform way. However, such a bound of convexity is 
difficult to obtain in general cases and may be even impossible as for the 
example studied in Section \ref{section-appli2}.

A typical example of application is given in Figure \ref{fig-foliation-1}: if 
$\Omega$ is a sphere and $\omega$ covers more than an hemisphere, then if $z$ 
is a global solution of a linear wave equation which vanishes in $\omega$ for 
all times, then $z\equiv 0$. Notice that, in this case, the family of surfaces 
is uniformly pseudo-convex and the unique continuation holds in fact in finite 
time even if $z$ is not a global in time solution.

\subsection{Unique continuation through pseudo-convex surfaces with boundary}

The case where the pseudo-convex surfaces $\Sigma_\lambda$ meet the boundary 
$\partial\Omega$ is more involved. Theorem \ref{th-hormander} has been 
generalized to this case by Tataru (see \cite{Tataru_uniq_analyticCPDE, Tataru, 
Tataru_uniq_analyticJMPA}). 
The boundary conditions are more difficult to describe geometrically, so we 
will 
only deal here with the case of flat geometry, that is $g(x,\xi)=|\xi|^2$, and  
the case of Dirichlet boundary condition. 

\begin{figure}[ht]

\begin{center}
\resizebox{6cm}{!}{\input{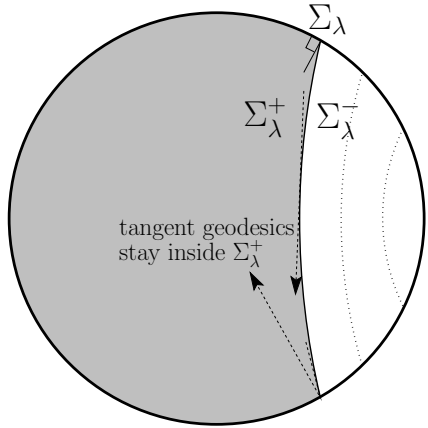}}
\end{center}
\caption{\it A oriented pseudo-convex foliation with boundary in the disk. The 
surfaces $\Sigma_\lambda$ forms a smooth family of curves inside a 
semidisk. The surfaces $\Sigma_\lambda^-$ are strictly convex and the angle 
formed by $\Sigma_\lambda$ and the boundary of the disk is less than $\pi/2$ on 
$\Sigma_\lambda^-$ side. The geodesics are straight lines bouncing at the 
boundary according to Newton's laws. The ones which are tangent 
to a surface $\Sigma_\lambda$ stay inside 
$\Sigma_\lambda^+$.\label{fig-foliation-2}}
\end{figure}

\begin{defi}\label{defi-foliation2}
A family of surfaces $(\Sigma_\lambda)_{\lambda\in[0,1)}$ is an \emph{oriented 
pseudo-convex foliation with boundary} in a flat manifold $\Omega$ if:
\begin{enumerate}[(i)]
\item the family of surfaces is smooth in the sense that it is locally 
described 
as level sets $\{x,\psi_\lambda(x)=0\}$ where 
$(x,\lambda)\mapsto\psi_\lambda(x)$ is a local smooth function with 
$\grad_x\psi_\lambda\neq 0$.
\item each surface is globally oriented in the sense that 
there exist disjoint sets $\Sigma_\lambda^\pm$ such 
that locally $\{x\in\Omega,\pm\psi_\lambda(x)>
0\}\subset \Sigma_\lambda^\pm$ and such that $\Omega=\Sigma_\lambda^- \cup 
\Sigma_\lambda \cup \Sigma_\lambda^+$.
\item for each $\lambda$, $(x,t)\mapsto \psi_\lambda(x)$ is pseudo-convex in 
the 
sense of definition  \ref{defi-pseudoconvex} as a function independent of $t$. 
Equivalently, $\Sigma_\lambda^-$ is locally strictly convex in a neighborhood 
of its boundary: the tangent space to $\Sigma_\lambda$ at $x_0$ is locally 
included in $\Sigma_\lambda^+$, $x_0$ excepted.
\item if a surface $\Sigma_\lambda$ meet $\partial\Omega$ at $x$, then 
$\partial_\nu \psi_\lambda(x)<0$. Equivalently, the angle formed by 
$\Sigma_\lambda$ and $\partial\Omega$ in the  region $\Sigma_\lambda^-$ is 
strictly less than $\pi/2$.
\end{enumerate}
\end{defi}
A typical example of such oriented pseudo-convex foliation with boundary is 
given in Figure \ref{fig-foliation-2}. Notice that the condition at the 
boundary is consistent with the one inside the domain. Indeed, the geodesics 
are straight lines which bounce at the boundary according to Newton's laws. 
Geometrically, we ask that any geodesic either crosses $\Sigma_\lambda$ in a 
transversal way, or stay locally inside $\Sigma_\lambda^+$.

By the same arguments as the ones in the proof of Theorem 
\ref{th-hormander-global} and using the result of Tataru, we obtain a global 
unique continuation result. 
\begin{theorem}\label{th-hormander-global-2}
Let $\Omega\subset\RR^d$ be a compact domain and let $\omega\subset\Omega$ be 
an open set. Assume that there exists an oriented pseudo-convex foliation 
$(\Sigma_\lambda)_{\lambda\in[0,1)}$ of $\Omega$ in the sense of Definition 
\ref{defi-foliation2}. Also assume that $\Sigma_0^+\subset 
\omega$ and $\omega\cup\left(\bigcup_{\lambda\in[0,1)} \Sigma_\lambda^+\right)$ 
covers $\Omega$ up to a set of zero measure. 

Let $b$, $c$ and $d$ be bounded  coefficients. Assume that $z$ is a global mild 
solution of
$$
\left\{\begin{array}{ll}
\partial^2_{tt} z=\Delta z + b(x,t)\partial_t z+c(x,t).\grad
z+d(x,t)z~~~~~&(x,t)\in\Omega\times \RR\\
z_{|\partial\Omega}(x,t)=0 & (x,t)\in\partial\Omega\times\RR\\
z_{|\omega}(x,t)=0&(x,t)\in\omega\times\RR
\end{array}
\right.
$$
Then $z\equiv 0$ everywhere.
\end{theorem}

A typical example of application is given in Figure \ref{fig-foliation-2}: if 
$\Omega$ is a disk and $\omega$ covers more than half of the boundary, then if 
$z$ is a global solution of a linear wave equation which vanishes in $\omega$ 
for all times, then $z\equiv 0$. 

\subsection{Proof of Theorem \ref{th2}}
Theorem \ref{th2} is then a direct consequence of the unique continuation 
results stated in this Section: Proposition \ref{prop-UCP-implique-th} and 
Theorems \ref{th-hormander-global} and \ref{th-hormander-global-2} imply 
Theorem \ref{th2}.

\section{Application 2: the peanut of rotation}\label{section-appli2}

We consider in this section the example of the peanut of rotation: a 
two-dimensional manifold where a central part is equivalent to the 
cylinder $\{x=(y,\theta)\in (-1,1)\times \SS \}$ endowed with the metric 
$g(y,\theta)={\rm d}y^2+\cosh^2(y){\rm d}\theta^2$ (see Figure \ref{fig-1}). 
The damping $\gamma$ is assumed to be positive, except in a part $x\in 
(-\ell,\ell)$ of the central part ($\ell\in (0,1)$). The 
decay of the linear damped wave semigroup has been established in \cite{CSVW} 
and \cite{Schenck}.
\begin{theorem}\label{th-CSVW}
{\bf Christianson, Schenck, Vasy \& Wunsch, 2014.}\\
In the setting of the peanut of rotation, there exist two positive constants 
$C$ and $\lambda$ such that the semigroup $e^{At}$ satisfies
$$\|e^{At} U_0\|_X \leq  C e^{-\lambda\sqrt{t}} \|U_0\|_{D(A)}~.$$
\end{theorem}

\begin{figure}[ht]
\begin{center}
\resizebox{8cm}{!}{\input{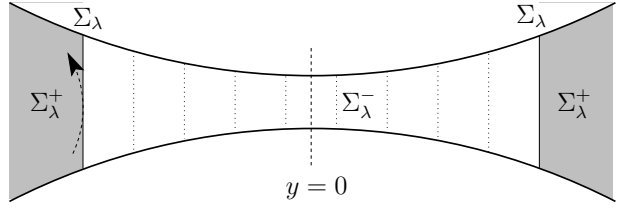}}
\end{center}
\caption{\it An oriented pseudo-convex foliation without boundary of the 
central part of the peanut. Each surface $\Sigma_\lambda$ consists in two 
vertical circles, the set $\Sigma_\lambda^-$ being the interior part surrounded 
by these circles. Since the central part of the peanut is 
negatively curved, the geodesics tangent to the vertical circles stay in the 
exterior domain $\Sigma_\lambda^+$. Notice that the set $\Sigma_0$ consists in 
both exterior disk and the disks $\Sigma_\lambda$ get closer to $x=0$ when 
$\lambda$ get closer to $1$. Thus, the central circle $x=0$ is not 
included in $\bigcup_{\lambda\in [0,1)}\Sigma_\lambda^+$ but this is not 
important since this set is of zero measure. 
\label{fig-foliation-peanut}}
\end{figure}

The decay rate of Theorem \ref{th-CSVW} obviously satisfies \eqref{hyp-decay-1} 
and \eqref{hyp-decay-2}. Thus, once the unique continuation property is 
obtained, Proposition \ref{prop-UCP-implique-th} yields the conclusion 
of Theorem \ref{th1} for the framework of the peanut of rotation. To obtain the 
unique continuation property, we will apply Theorem \ref{th-hormander-global} 
with the family of pseudo-convex surfaces shown in Figure 
\ref{fig-foliation-peanut}.

\noindent Applying Theorem \ref{th2} and the ideas of Proposition \ref{prop1}, 
we obtain the following result.
\begin{theorem}\label{th-appli2}
Consider the damped wave equation \eqref{eq} in the framework of the peanut of 
rotation introduced above. Let $\alpha>0$ and $f$ satisfying 
\eqref{hyp-f-estim} and \eqref{hyp-f-sign}. Then, any solution $u$ of 
\eqref{eq} satisfies 
$$\|(u,\partial_t u)(t)\|_{H^1_0\times L^2}  ~\xrightarrow[~~t\longrightarrow 
+\infty~~]{}~0~.$$

Moreover, for any $R$ and $\sigma\in (0,1]$, there exists 
$C_{R,\sigma}$ such that, for any solution $u$ with $U_0\in 
X^\sigma$, 
$$ \|(u_0,u_1)\|_{H^{1+\sigma}\times H^\sigma}\leq 
R~~\Longrightarrow ~~ \|(u,\partial_t u)(t)\|_{H^1_0\times L^2} \leq 
C_{R,\sigma}e^{-\sigma\tilde \lambda\sqrt{t}}$$
where $\tilde\lambda$ is the linearized rate given in Lemma \ref{lemma-appli2}.
\end{theorem}
\begin{demo}
Let us first formally check that the family of disks introduced in Figure 
\ref{fig-foliation-peanut} is a suitable pseudo-convex foliation without 
boundary. We use the cylindrical coordinates $x=(y,\theta)\in (-1,1)\times\SS$ 
with associated tangent variables $\xi=(\zeta,\Theta)$.
By symmetry, we only consider the right-hand-side circles which are 
defined by $\psi_\lambda(x)=0$ with 
$\psi_\lambda(x)=y-(1-\lambda)$. The circle $\Sigma_0$ corresponding to 
$y=1$ is in the interior of the region $\omega$ where the 
damping is positive. When $\lambda$ get closer to $1$, the circle 
$\Sigma_\lambda$ get closer to $y=0$. We are obviously in the setting of 
Definition \ref{defi-foliation1} and thus of Theorem \ref{th2}, except maybe 
for the assumption of strong pseudo-convexity. We already give a geometrical 
insight of this assumption, but let us check it formally. 

The local metric is given by $g(y,\theta)={\rm d}y^2+\cosh^2(y){\rm 
d}\theta^2$. The Laplace-Beltrami operator is thus given by
$$\Delta_g = \partial_{yy}^2 + 2\tanh(y)\partial_y + \frac 1 
{\cosh^2(y)}\partial^2_{\theta\theta}~.$$
The principal part of the wave operator is then
$$p(y,\theta,t,\zeta,\Theta,
\tau)=|\zeta|^2+\frac 1{\cosh^2(y)}|\Theta|^2-|\tau|^2~.$$
Thus $H_p(\psi_\lambda)=2\zeta$ and 
$$H_p^2(\psi_\lambda)=4 \frac{\sinh(y)}{\cosh(y)^3} |\Theta|^2~.$$
The pseudo-convexity condition is then checked. Indeed, if 
$H_p(\psi_\lambda)=0$ then $\zeta=0$ and since $\xi=(\zeta,\Theta)$ must be 
non-zero, we must have $\Theta\neq 0$. As $\psi_\lambda(y,\theta)=0$ with 
$\lambda\in[0,1)$, we have $y>0$ and thus $H_p^2(\psi_\lambda)>0$. Looking 
carefully to the computations, one notes that, in fact, we only need that the 
radius $\cosh(y)$ of the cylindrical part is increasing for $y>0$ and 
decreasing for $y<0$ to obtain the unique continuation property.

The above arguments show the stabilization of the semilinear damped wave 
equation. Moreover, Proposition \ref{prop-asymp-compact-2} shows the uniform 
convergence to $0$ in $X^{\sigma'}$ for initial data in a more regular space 
$X^\sigma$ ($\sigma>\sigma'$). 

To obtain the decay estimate of Theorem 
\ref{th-appli2}, we argue as in the proof of the second case of 
Proposition \ref{prop1}. However, we claim that we can avoid the loss in the 
power by following the remark below Lemma \ref{lemma-appli2}. Indeed, Theorem 
5.1 of \cite{CSVW} implies in our framework that, for large $\mu$,
$$\|(-\Delta+\mu^2-i\mu\gamma)^{-1}\|_{\Lc(L^2)}\,\leq\, 
C\,\frac{\ln|\mu|}{|\mu|}~.$$
We argue as in the proof of Lemma \ref{lemma-appli2}. Applying Propositions 
\ref{prop-resolvent}, \ref{prop-estimP} and \ref{prop-appendix-linearized}, we 
obtain 
$$\|(\tilde A-i\mu)^{-1}\|_{\Lc(L^2)} \leq C \, {|\ln \mu|}~.$$
and by Theorem \ref{th-BD}, we get the conclusion of Lemma \ref{lemma-appli2} 
in the form 
$$\|e^{\tilde A t}U_0\|_X \leq Ce^{-b\sqrt{t}} \|U_0\|_{D(A)}~.$$
In other words, we avoid an additional logarithmic loss by directly dealing 
with the estimates of \cite{CSVW} instead of using the back and forth 
implications of Theorems \ref{th-BD-1} and \ref{th-BD}.

It is then sufficient to follow the proof of the exponential case of 
Proposition \ref{prop1} with the exponent $\gamma=1/2$. 
\end{demo}

\section{Decay estimate in the disk with holes}\label{section-decay-disk}
In the previous examples of application, the decay of the semigroups was 
explicitly written in previous papers. In the case of a disk with several 
holes, we are not aware of a paper where an explicit decay is written. The 
corresponding scattering problem has been studied by Ikawa in 
\cite{Ikawa1,Ikawa2}. Many further studies have been published. In this 
article, we will use an estimate and a ``black box argument'' introduced by 
Burq and Zworski in \cite{Burq-Zworski}. Combining them with the results in 
Appendix, we obtain the following decay.

\begin{theorem}\label{th-decay-disque-trous}
Let $\Oc\subset\RR^d$ be a smooth bounded open set. For $i=1\ldots p$, let 
$O_i\subset \Oc$ be smooth strictly convex obstacles satisfying:
\begin{enumerate}[(a)]
\item the obstacles are disjoint: $O_i\cap O_j=\emptyset$ for $i\neq j$,
\item the convex hull {\it convhull}$(\cup_i O_i)$ of the obstacles is 
contained in $\Oc$,
\item no obstacle is in the convex hull of two others, that is that {\it 
convhull}$(O_i\cup O_j)\cap O_k =\emptyset$ for $i,j,k$ different,
\item if there are three or more obstacles ($p>2$), set $\kappa$ the infimum of 
the principal curvatures of the boundaries $\partial O_i$ of the obstacles and 
$L$ the minimal distance between two obstacles, and assume that $\kappa L>p$.  
\end{enumerate}
Let $O=\cup_i O_i$, let $\Omega=\Oc\setminus O$ be the domain with 
holes and let $\gamma\geq 0$ be a damping which is strictly positive in a 
neighborhood of the exterior boundary $\partial \Oc$. Then the semigroup 
$e^{At}$ of the linear damped wave equation on $\Omega$ satisfies
$$\|e^{At}U_0\|_X \leq C e^{-\lambda t^{1/3}} \|U_0\|_{D(A)}$$
with $C$ and $\lambda$ two positive constants.
\end{theorem}

\begin{figure}[ht]
\begin{center}
\resizebox{6cm}{!}{\input{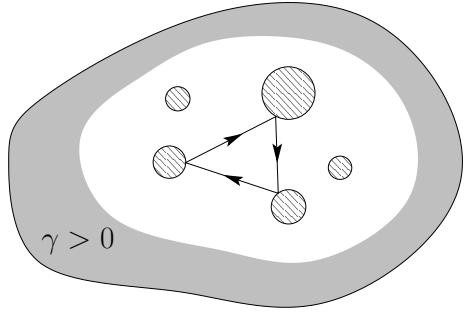}}
\end{center}
\caption{\it A domain satisfying Hypotheses (a)-(d) of 
Theorem \ref{th-decay-disque-trous}. Notice the presence of 
periodic geodesics bouncing on the obstacles and never meeting the support of 
the damping $\gamma$. The geometric control condition thus fails and we cannot 
hope a uniform decay of the semigroup in $X$. \label{fig-disque-trous}}
\end{figure}

A typical domain consists in a smooth domain with several small holes as in 
Figure \ref{fig-disque-trous}. Typically, if $O_i$ are balls of center $c_i$ 
and radius $r_i$ and if there is no triplet of aligned centers, then Assumption 
(d) holds for $r_i$ small enough since $\kappa$ becomes large whereas $L$ stay 
bounded.\\

\noindent {\bf Remarks:}
\begin{itemize}
\item We do not claim that the decay rate $e^{-t^{1/3}}$ is optimal. In fact, 
our proof uses rough arguments leading to logarithmic losses. We 
strongly believe that the right decay rate is $e^{-t^{1/2}}$. A strategy of 
proof may be to follow the arguments of Datchev and Vasy \cite{Vasy}, adding 
the presence of a boundary. However, this improvement is not central and would use 
techniques too far from the spirit of this paper.
\item It is also certainly possible to relax the assumption about $\gamma$ to 
the following weaker assumption. There exists a neighborhood $K$ of $O=\cup_i 
O_i$ such that any geodesic ray starting in $\Omega\setminus K$ reach a point 
where $\gamma\geq \varepsilon >0$ before meeting $K$, for at least the 
backward or the forward flow.
\end{itemize}

Based on \cite{Ikawa2}, the following estimate appears in \cite{Burq-Zworski}
\begin{prop}\label{prop-BZ}
{\bf Ikawa's black box (Section 6.2 of \cite{Burq-Zworski})}\\
Let $O_i$ be obstacles in $\RR^d$ satisfying the Assumptions (a),(c) and (d) of 
Theorem \ref{th-decay-disque-trous}. 
Let $R_{bb}(\mu)$ be the outgoing resolvent of the Laplacian operator outside 
the obstacles, that is the meromorphic continuation of 
$(\Delta_{\RR^d\setminus O}+\mu^2)^{-1}$ from Im$(\mu)>0$, where 
$\Delta_{\RR^d\setminus O}$ is the Laplacian operator on the 
exterior domain ${\RR^d\setminus 
O}=\RR^d\setminus (\cup_i O_i)$ with Dirichlet boundary condition.

Then, for any cut-off function $\chi\in\Cc^\infty_c(\RR^d)$, we have 
$$\|\chi R_{bb}(\mu) \chi\|_{\Lc(L^2({\RR^d\setminus 
O}))} \leq C \frac {\ln (1+|\mu|)}{1+|\mu|}~. 
$$
\end{prop}

Using this black box in the same spirit of \cite{Burq-Zworski}, we obtain the 
following observation estimate.
\begin{lemma}\label{lemma-BZ}
Assume that the assumptions of Theorem \ref{th-decay-disque-trous} hold. Then 
there exists $C>0$ such that
$$\forall \mu\in\RR~,~~\forall u\in D(\Delta)~,~~\|u\|_{L^2(\Omega)} ~ \leq ~  
C\frac {\ln |\mu|}{|\mu|} \|(\Delta+\mu^2)u\|_{L^2(\Omega)} + C{\ln |\mu|} 
\|\sqrt{\gamma} u\|_{L^2(\Omega)}$$
where $\Delta$ is the Laplacian operator on the bounded 
domain with holes $\Omega=\Oc\setminus O$ with Dirichlet boundary conditions.
\end{lemma}
\begin{demo}
By compactness, there exists $\varepsilon>0$ such that $\gamma\geq \varepsilon 
>0$ in a neighborhood of the exterior boundary $\partial \Oc$. Let 
$\chi\in\Cc^\infty_0(\RR^d)$ be a smooth cut-off function equal to $1$ in a 
neighborhood of {\it convhull}$(\cup_i O_i)$, equal to $0$ outside $\Oc$ and 
such that $1-\chi$ is supported where $\gamma\geq \varepsilon >0$. This is 
possible by Assumption (b). We have immediately
\begin{equation}\label{eq-lemma-BZ-1}
\forall u\in L^2(\Omega)~,~~\|(1-\chi) u\|_{L^2(\Omega)} \leq \frac 
1{\sqrt{\varepsilon}} \|\sqrt{\gamma} u\|_{L^2(\Omega)}~.
\end{equation}
Let $\chi_0$ be a smooth cut-off function supported in $\Omega$ so that 
$\chi_0\equiv 1$ in a neighborhood of the support of $\chi$ and let $\chi_1$ be 
another smooth cut-off function supported in $\Omega$ so that 
$\chi_1\equiv 1$ in a neighborhood of the support of $\chi_0$. For all $u\in 
H^2(\Omega)\cap H^1_0(\Omega)$, we extend $\chi_i u$ as a function in 
$H^2(\RR^d\setminus O)\cap H^1_0(\Omega)$. Regarding $\chi_i u$, applying 
$\Delta_{\RR^d\setminus O}$ or $\Delta$ gives the same result. We can thus 
apply the ``black box'' estimate of Proposition \ref{prop-BZ} as follows.
\begin{align*}
\|\chi u\|_{L^2(\Omega)} & = \|\chi \chi_0^2 u\|_{L^2(\Omega)} =  \|\chi\chi_0 
R_{bb}(\mu)(\Delta+\mu^2)\chi_0 u\|_{L^2(\Omega)}\\
&\leq \|\chi\chi_0 
R_{bb}(\mu)\chi_0(\Delta+\mu^2) u\|_{L^2(\Omega)} + \|\chi\chi_0 
 R_{bb}(\mu) [\Delta,\chi_0] u\|_{L^2(\Omega)}\\
&\leq \|\chi_0 R_{bb}(\mu)\chi_0\|_{\Lc(L^2(\RR^d\setminus O))} 
\|(\Delta+\mu^2)u\| + \|\chi\chi_0 
\chi_1 R_{bb}(\mu)\chi_1 [(\Delta+\mu^2),\chi_0] u\|\\
&\leq \|\chi_0 R_{bb}(\mu)\chi_0\|_{\Lc(L^2(\RR^d\setminus O))} 
\|(\Delta+\mu^2) u\| + \|\chi_1 
R_{bb}(\mu)\chi_1\|_{\Lc(L^2(\RR^d\setminus O))}
\|[\Delta,\chi_0] u\|\\
&\leq C \frac {\ln (1+|\mu|)}{1+|\mu|}\left( \|(\Delta+\mu^2) u\|_{L^2(\Omega)} 
+ \|u\|_{H^1(\text{supp}(\grad \chi_0))}\right)~.
\end{align*}
By interpolation and elliptic regularity, we have
\begin{align*}
\|u\|^2_{H^1(\text{supp}(\grad \chi_0))} & \leq C \|u\|_{L^2(\text{supp}(\grad 
\chi_0))} \|u\|_{H^2(\text{supp}(\grad \chi_0))}\\
&\leq C \|u\|_{L^2(\text{supp}(\grad \chi_0))}\left(\|\Delta 
u\|_{L^2(\text{supp}(\grad \chi_0))} + \|u\|_{L^2(\text{supp}(\grad 
\chi_0))}\right)\\
&\leq C \|u\|_{L^2(\text{supp}(\grad 
\chi_0))}\left(\|(\Delta+\mu^2)u\|_{L^2(\text{supp}(\grad \chi_0))} + 
(1+\mu^2)\|u\|_{L^2(\text{supp}(\grad \chi_0))}\right)\\
& \leq C \left(\frac 1{2(1+\mu^2)}\|(\Delta+\mu^2)u\|_{L^2(\text{supp}(\grad 
\chi_0))} + 
\frac 32(1+\mu^2)\|u\|_{L^2(\text{supp}(\grad \chi_0))}\right)
\end{align*}
Since the support of $\grad\chi_0$ is included in the place where $\gamma\geq 
\varepsilon >0$, both previous estimates yield
$$\|\chi u\|_{L^2(\Omega)} \leq C  \frac {\ln (1+|\mu|)}{1+|\mu|}\left( 
\|(\Delta+\mu^2) u\|_{L^2(\Omega)} 
+ (1+|\mu|)\|\sqrt \gamma u\|_{L^2(\Omega)}\right)~.$$
With \eqref{eq-lemma-BZ-1}, this concludes the proof.
\end{demo}

{ \noindent \emph{\textbf{Proof of Theorem 
\ref{th-decay-disque-trous}:}}}
Applying Propositions \ref{prop-1-appendix} and \ref{prop-estimP} in 
Appendix, the observability estimate of Lemma \ref{lemma-BZ} implies 
that there exists $C>0$ such that  
$$\forall \mu\in\RR~,~~ \|(A-i\mu)\|_{\Lc(X)} ~\leq~ C \ln^2(2+|\mu|)~.$$
Then, we apply Theorem \ref{th-BD} of Batty and Duyckaerts stated in Appendix 
to obtain the decay with rate $e^{-\lambda t^{1/3}}$. Notice that Proposition 
\ref{prop-1-appendix} and Theorem \ref{th-BD} contain some losses transforming 
the rate $\ln |\mu|$ of Lemma \ref{lemma-BZ} into first $\ln^2 |\mu|$ and then 
$\ln^3 |\mu|$. This is responsible of the power $1/3$ in the decay rate. 
{\hfill$\square$\\}

\section{Application 3: the disk with two holes}\label{section-appli3}

\begin{figure}[ht]
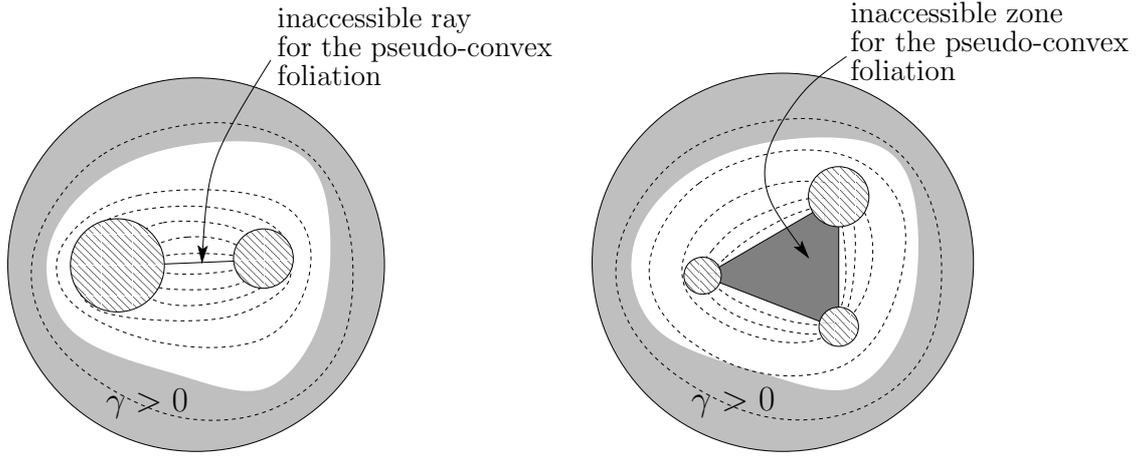

\begin{center}
\resizebox{7.5cm}{!}{\input{disque-trous-3.pstex_t}}$~$\resizebox{7.5cm}{!}{
\input { disque-trous-4.pstex_t}}
\end{center}
\caption{\it Two examples of disks with holes and associated attempts to draw a 
suitable pseudo-convex foliation covering the whole domain. 
Left, the disk with two holes may be covered by a pseudo-convex foliation 
starting in a neighborhood of the exterior boundary, except for a single line, 
which is of measure zero. Right, an attempt to cover a disk with three holes 
with a pseudo-convex foliation. We easily notice that the central zone cannot 
be covered and Theorem \ref{th-hormander-global-2} cannot be applied in this 
case.}\label{fig-appli3}
\end{figure}

In the previous section, we have obtained a sufficiently fast decay rate for 
the semigroup of the damped wave equation in a disk with several holes as in 
Figure \ref{fig-1}. If we prove the unique continuation property of 
Proposition \ref{prop-UCP-implique-th} in this situation, then we would obtain 
the desired stabilization. To obtain the unique continuation property, we would 
like to use Theorem \ref{th-hormander-global-2}, that is to exhibit an oriented 
pseudo-convex foliation $(\Sigma_\lambda)_{\lambda\in[0,1)}$ with $\Sigma^+_0$ 
included in a neighborhood of the boundary and $\cup_\lambda \Sigma^+\infty$ 
covering almost all $\Omega$. This is possible in the case where there is at 
most two holes in the disk and impossible if there are more holes, as shown in 
Figure \ref{fig-appli3}.

Thus, in the case where there is only two holes, Theorem 
\ref{th-hormander-global-2} enables to use Proposition 
\ref{prop-UCP-implique-th} and to obtain the conclusions of Theorem \ref{th1}. 
Moreover, notice that in this case, there is no technical assumptions, neither 
(c) nor (d), in Theorem \ref{th-decay-disque-trous}. We thus obtain the 
following result, as an application of Theorem \ref{th2}, 
Proposition \ref{prop1} and mutatis mutandis the same use of the remark 
below Lemma \ref{lemma-appli2} as in the proof Theorem \ref{th-appli2}.  
The proof is left to the reader.
\begin{theorem}\label{th-appli3}
Consider the damped wave equation \eqref{eq} in the framework of a disk with 
two convex holes and assume that the damping $\gamma$ is strictly positive in a 
neighborhood of the exterior boundary. Let $\alpha>0$ and $f$ satisfying 
\eqref{hyp-f-estim} and \eqref{hyp-f-sign}. Then, any solution $u$ of 
\eqref{eq} satisfies 
$$\|(u,\partial_t u)(t)\|_{H^1_0\times L^2}  ~\xrightarrow[~~t\longrightarrow 
+\infty~~]{}~0~.$$

Moreover, there exists $\tilde\lambda$ such that, for any $R$ and $\sigma\in 
(0,1]$, there exists $C_{R,\sigma}$ such that, for any solution $u$ with 
$U_0\in X^\sigma$, 
$$ \|(u_0,u_1)\|_{H^{1+\sigma}\times H^\sigma}\leq 
R~~\Longrightarrow ~~ \|(u,\partial_t u)(t)\|_{H^1_0\times L^2} \leq 
C_{R,\sigma}e^{-\sigma\tilde \lambda{t}^{1/3}}~.$$
\end{theorem}

\section{Analytic regularization and proof of 
Theorem \ref{th1}}\label{section-reg-analytic}

In Figure \ref{fig-appli3}, we have seen that, in some situations, the unique 
continuation property of Lerner-Robbiano-H\"ormander stated 
in Theorem \ref{th-hormander-global-2} is not useful. To deal with these 
situations, we need another unique continuation property: Theorem \ref{th-RZ} of 
Robbiano-Zuily-H\"ormander to \eqref{eq-UCP}. This is only possible if the 
coefficients of \eqref{eq-UCP} are analytic with respect to the time $t$. Thus, 
we need to choose $f'_u$ analytic with respect to $u$ and to prove 
that the global solution $w$ appearing in \eqref{eq-UCP} is analytic in time. 
This is the basic idea leading to Theorem \ref{th1}.

However, even if $f$ is analytic, the damped wave 
equation does not regularize its solutions. To overcome this problem,  
we use the following fact known since the work of Hale 
and Raugel in \cite{Hale-Raugel}: the globally bounded solutions of the 
damped wave equation are as smooth as the non-linearity $f$. This asymptotic 
regularization property is linked to the asymptotic smoothness or 
compactness property (see Section \ref{section-asymp-compact}). The idea that 
this asymptotic smoothing of the damped wave equation may be used to apply 
analytic unique continuation theorems originates from the work of Hale and 
Raugel, even if they did not publish this idea. The first 
published occurrence appears in \cite{RJ} (see also \cite{RC}).

The article \cite{Hale-Raugel} contains several abstract theorems. They apply 
for linear semigroup with uniform decay, that is $\|e^{At}\|_{\Lc(X)}\leq M 
e^{-\lambda t}$. The purpose of the present article is to study cases where 
this uniform decay fails, so \cite{Hale-Raugel} does not directly apply: we need 
to extend its results in our case where the semigroup has a weaker decay. 
Extending these results in the most general framework will lead to heavy 
notations and assumptions. That is why, we only consider here damped wave 
equations in a simple setting and in low dimension.

\subsection{Analytic regularization of global bounded solutions}
Let $d=2$ or $3$ and let $\Omega$ be a smooth manifold of dimension $d$ 
with or without boundary and such that $\overline\Omega$ is compact. Let 
$\gamma\in \Cc^1(\Omega,\RR)$ be a nonnegative damping, let $\Delta$ be the 
Laplacian operator with Dirichlet boundary condition and let 
$f\in\Cc^\infty(\overline\Omega\times\RR,\RR)$ be a smooth nonlinearity. We 
assume that $f$ is of polynomial type in the sense that there exist 
$p>0$ and $C>0$ such that \eqref{hyp-f-estim} holds.
We consider global solutions of the damped wave equation 
\begin{equation}\label{eq-analytic}
\left\{\begin{array}{ll}
\partial_{tt}^2 u(x,t) + \gamma(x)\partial_t u(x,t) = \Delta u(x,t) - \alpha 
u(x,t) -  f(x,u(x,t))~~~&(x,t)\in\Omega\times \RR\\
u_{|\partial\Omega}(x,t)=0&(x,t)\in\partial\Omega\times \RR
\end{array}\right.
\end{equation}
We use the notations of Section \ref{section-basic}. In particular, for $\sigma\in[0,1]$, 
we set $X^{\sigma}=(H^{1+\sigma}(\Omega)\cap H^1_0(\Omega))\times
H^{\sigma}_0(\Omega)$, $X^0=X=H^1_0(\Omega)\times L^2(\Omega)$ 
and $X^1=D(A)$. Also notice that $H^{1+\sigma}(\Omega)$ is a subspace of 
$\Cc^0(\overline\Omega)$ for $d=2$ and $\sigma>0$ or for $d=3$ and 
$\sigma>1/2$. 

The purpose of this section is to prove the following result
\begin{theorem}\label{th-analytic}
Assume that the above setting holds, in particular assume that $\gamma\in \Cc^1(\Omega,\RR)$. 
Let $U(t)\in\Cc^0(\RR,X)$ be a mild solution of \eqref{eq-analytic} and assume moreover that 
\begin{enumerate}[(i)]
\item There exists $\sigma_0\in (0,1)$ such that $U(t)$ is defined for all 
$t\in\RR$ and uniformly bounded in $X^{\sigma_0}$, that is that there 
exists $C>0$ such that
$$\forall t\in\RR~,~~ \|U(t)\|_{X^{\sigma_0}}\leq C~.$$
If $d=3$, assume in addition that $\sigma_0>1/2$.
\item The linear semigroup satisfies the decay estimate 
\begin{equation}\label{hyp-decay-analytic}
\forall U_0\in D(A)~,~\|e^{At}U_0\|_X ~\leq ~\frac M 
{t^{\beta}}\,\|U_0\|_{D(A)}
\end{equation}
with $\beta>\max(2p,2/(1-\sigma_0))$, where $p$ is the polynomial growth of $f$ in 
\eqref{hyp-f-estim}.
\item The function $u\in\RR \mapsto f(x,u)\in\RR$ is analytic with respect to 
$u$, 
\end{enumerate}
Then the mapping $t\in\RR \mapsto U(t)\in X^{\sigma_0}$ is analytic with 
respect to $t$.
\end{theorem}

We expect that the condition $\beta>2p$ in Assumption (ii) is not optimal: at 
least $\beta>p$ should be sufficient if we get rid of the losses in the too 
general proofs of the auxiliary results in appendix and $\beta>2/(1-\sigma_0)$ 
could be omitted by assuming more regularity on $\gamma$. 
Since our concrete applications have a linear decay of the type  
$\Oc(e^{-t^{\beta}})$, we let these probable improvements for later study. 

\vspace{3mm}

The remaining part of this section is devoted to the proof of Theorem 
\ref{th-analytic}. 

\vspace{3mm}

{\noindent \bf $\bullet$ Step 1: a trick to satisfy the boundary 
condition}\\[1mm]
We would like that $u\mapsto f(\cdot,u)$ maps $H^{1+\sigma}(\Omega)\cap 
H^1_0(\Omega)$ into $H^\sigma_0(\Omega)$. For $d=2$ or $d=3$ and $\sigma>1/2$, 
the regularity part is trivial since $H^{1+\sigma}(\Omega)\subset \Cc^0(\Omega)$ 
is a normed algebra. However, the boundary condition is not necessarily 
fulfilled since $f(x,0)$ may be different from $0$. Let us describe here a 
trick introduced in \cite{Hale-Raugel} to deal with this problem. Let $e$ be 
the (time-independent) solution of $\Delta e -\alpha e = f(x,0)$ in 
$H^2(\Omega)\cap H^1_0(\Omega)$. We notice that $\tilde u=u-e$ solves
$$\partial^2_{tt} \tilde u + \gamma(x)\partial_t \tilde u = \Delta \tilde u - 
\alpha \tilde u - f(x,\tilde u + e(x)) + f(x,0)~.$$
If we set 
$$\tilde f(x,\tilde u)= f(x,\tilde u + e(x)) - f(x,0)$$
we obtain a function $\tilde f$ as smooth as $f$ which is also analytic with 
respect to $\tilde u$. Moreover, $\tilde f(x,0)=0$ for $x\in\partial\Omega$, 
which shows that $\tilde u\mapsto \tilde f(\cdot,\tilde u)$ maps 
$H^{1+\sigma}(\Omega)\cap H^1_0(\Omega)$ into $H^\sigma_0(\Omega)$, including 
the boundary condition. If we prove Theorem \ref{th-analytic} for $\tilde u$ 
and $\tilde f$, it clearly yields Theorem \ref{th-analytic} for $u$ and $f$. 
Thus, we may assume that $f(x,0)=0$ at the boundary and we will forget the 
tilde sign to lighten the notations in what follows. In particular, $F=(0,f)$ 
maps $X^\sigma$ into $X^{1+\sigma}$.

\vspace{3mm}

{\noindent \bf $\bullet$ Step 2: the decay of the high-frequencies 
semigroups}\\[1mm]
From now on, we also use the notations of the appendices, Sections 
\ref{section-appendix1} and \ref{section-appendix2}. In particular, we set
$L=-\Delta+\alpha$. Since $L$ is self-adjoint, positive and with compact 
resolvent, there exists an 
orthonormal basis $(\phi_k)_{k\geq 0}$ of eigenfunctions of $L$ corresponding to 
the eigenvalues $(\lambda_k)_{k\geq 0}$. As in Section \ref{section-appendix2}, 
we introduce the high-frequencies truncations $Q_n$, that are the projectors on 
the space span$\{\phi_k,k\geq n\}$ and we set $\Qc_n=(Q_n,Q_n)$ on $X$.

The proof of Theorem \ref{th-analytic} is based on the following splitting.
We introduce $P_n=Id-Q_n$ and $\Pc_n=Id-\Qc_n$, which are low-frequencies
projections with finite rank. We consider the splitting
\begin{equation}\label{eq-proj2}
U\,=\,\Pc_n U + \Qc_n U\,:=\,V+W~.
\end{equation}
Then \eqref{eq-analytic} writes
\begin{equation}\label{eq-proj}
\left\{
\begin{array}{l}\partial_t V = (\Pc_nA\Pc_n) V +
(\Pc_nA\Qc_n)W+\Pc_n F(V+W)\\
\partial_t W = (\Qc_nA\Qc_n) W +
(\Qc_nA\Pc_n)V+\Qc_n F(V+W)\\
\end{array}\right.
\end{equation}

As a consequence of the results in the appendices, we have the following 
decay estimates.
\begin{prop}\label{prop-ana-1}
Assume that Hypothesis (ii) of Theorem \ref{th-analytic} holds. Then, for all 
$\sigma\in[0,1)$, $\nu\in (0,1-\sigma]$ and $\varepsilon>0$, there exists 
$C>0$ such that, for all $n\in\NN$, 
\begin{equation}\label{eq-th-decay-2}
\forall t\geq 0~,~\forall U\in \Qc_n
X^{\sigma+\nu}~~\left\|e^{\Qc_nA\Qc_n
t} U\right\|_{X^\sigma}~\leq ~\frac C {t^{\nu 
\beta/2-\varepsilon}} \, \|U\|_{X^{\sigma+\nu}}~.
\end{equation}
\end{prop}
\begin{demo}
Using the arguments of Proposition \ref{prop-decay-interpole}, it is sufficient 
to obtain the decay for $\sigma=0$ and $\nu=1$, that is the decay 
estimate from $D(A)$ into $X$. 

In the appendices, we recall the result of Borichev and Tomilov in 
\cite{Bor-Tom} (see Theorem 
\ref{th-BT}). In the context of the damped wave equation, we may also consider 
Proposition 2.4 of \cite{Ana-Leau} by Anantharaman and L\'eautaud. We obtain 
that Hypothesis (ii) of Theorem \ref{th-analytic} implies the estimate
\begin{equation}
\|u\|_{L^2}\leq C \left( |\mu|^{1/\beta-1} \|P_B(\mu)u\|_{L^2} + 
|\mu|^{1/2\beta} \|\sqrt{B} u\| \right) 
\end{equation}
with $B$ is the multiplication by $\gamma$ and $P_B(\mu)=-L-i\mu B+\mu^2 Id$. Note that the unique 
continuation assumed in \cite[Proposition 2.4]{Ana-Leau} is satisfied because Hypothesis 
(ii) is satisfied. Thus, we obtain that 
\begin{equation}\label{eq-demo-analytic-1}
\|u\|_{L^2}\leq C \left( |\mu|^{1/\beta-1} \|P(\mu)u\|_{L^2} + 
\right(|\mu|^{1/\beta} \|\sqrt{B}\|_{\Lc(L^2)} + |\mu|^{1/2\beta}\left) 
\|\sqrt{B} u\| \right) 
\end{equation}
with $P(\mu)=-L+\mu^2 Id$.
Then, Propositions \ref{prop-1-appendix-high-frequencies} and 
\ref{prop-2-appendix} show that
$$\forall n\in\NN~,~~\forall \mu\in\RR~,~~\left\|(\Qc_n A \Qc_n - 
i\mu)^{-1}\right\|_{\Lc(\Qc_n X)} ~\leq ~ \frac K{|\mu|^b}$$
with $b=2/\beta$. At this point, we may use Theorem \ref{th-BT} in appendix to 
obtain the decay of the linear semigroup. However, this result of 
\cite{Bor-Tom} (as the one of \cite{Ana-Leau}) is not stated with explicit 
constants and we need to be sure that these constants are uniform in $n$ (even 
if this is surely the case). Thus, we accept here a small 
loss (harmless for the proof of Theorem \ref{th-analytic}) and we use the 
explicit statement of Theorem \ref{th-BD} to obtain that, for any 
$\varepsilon>0$, there exists a constant $C$ such that
$$\left\|e^{\Qc_nA\Qc_n t} U\right\|_{X}~\leq ~\frac C 
{t^{\beta/2-\varepsilon}} 
\, \|U\|_{D(A)}$$
which concludes the proof.
\end{demo}

\vspace{3mm}

{\noindent \bf $\bullet$ Step 3: the finite determining modes}\\[1mm]
In this step, we follow the arguments of \cite{Hale-Raugel} with the main modifications 
coming from the weaker decay of the linear semigroup. 
We consider the complex setting, that is that the functions in
$X$ are complex valued. We recall that a function $\Psi$ between two complex
Banach spaces $Y$ and $Z$ is said to be holomorphic if its Fr\'echet derivative
exists for any $y\in Y$. We introduce the notation
$$\Bc_{M,\delta}(Y) = \{U(t)\in \Cc^0(\RR,Y)~~|~~\forall
t\in\RR,~\|\Re(U(t))\|_Y\leq M\text{ and }\|\Im(U(t)\|_Y \leq \delta \}~.$$
The space $\Bc_{M,\delta}(Y)$ is naturally endowed with the $L^\infty(Y)$-norm.
We assumed that $\sigma_0>0$ for $d=2$ and $\sigma_0>1/2$ for $d=3$ as in 
Hypothesis (iii) of Theorem \ref{th-analytic}, so that 
$H^{1+\sigma_0}(\Omega)\subset \Cc^0(\Omega)$.

We will use the holomorphic extension of $f$ in a technical setting stated in the following
lemma. Except this particular setting, the result is a straightforward consequence of the analyticity
of $f$ and we omit the proof.
\begin{lemma}
Assume $u\in\RR \mapsto f(x,u)\in\RR$ is analytic with respect to 
$u$. Denote $\kappa$ the injection constant $\|\cdot\|_{L^\infty}\leq \kappa
\|\cdot\|_{H^{1+\sigma_0}}$.

Let $M$ be given and $M'_0\geq 0$. Then, there exists
$M'\geq M'_0$, as well as two small positive constants
$\delta$ and $\delta'$ such that the following holds. The function 
$z\in\RR \mapsto f(\cdot,z)\in \Cc^1(\overline\Omega,\RR)$ has a
holomorphic extension in $\{z\in\CC,~|\Re(z)|\leq \kappa(M+M'+\delta)$ and
$|\Im(z)|\leq \kappa(\delta+\delta')\}$.
\end{lemma}

We apply the above lemma to obtain the following result.
\begin{prop}\label{prop-W}
Let $M$ and $M'_0\geq 0$ be given and let $n_0\in\NN$. Let $M'$, 
$\delta$ and $\delta'$ the constants given by the previous lemma. Then, there exist
$n\geq n_0$ so that for any function $V(t)$ in the
complex set $\Bc_{M+\delta,\delta}(\Pc_nX^{\sigma_0})$, there exists a 
unique bounded solution $W$ in $\Bc_{M',\delta'}(\Qc_nX^{\sigma_0})$ of
\begin{equation}\label{eq-prop-W}
\partial_t W = (\Qc_nA\Qc_n) W + (\Qc_nA\Pc_n) V+\Qc_n F(V+W)~.
\end{equation}
In addition, the mapping $V\mapsto W(V)$ is lipschitzian and
holomorphic.
\end{prop}
\begin{demo}
Assume that $W$ solution of \eqref{eq-prop-W} exists and is bounded
in $\Cc^0(\RR,\Qc_nX^{\sigma_0})$. Then, 
\begin{align*} 
W(t)=&e^{\Qc_n A\Qc_n (t-t_0)}W(t_0) \\
&~~+~\int_{0}^{t-t_0} e^{\Qc_n A\Qc_n s}\left((\Qc_nA\Pc_n)
V(t-s)+\Qc_n F(V+W)(t-s)\right)\,{\rm d}s.
\end{align*}
Using Proposition \ref{prop-ana-1}, when $t_0$ goes to $-\infty$, we get
$$W(t)= \int_{0}^\infty e^{\Qc_n A\Qc_n s}\left((\Qc_nA\Pc_n)
V(t-s)+\Qc_n F( V+W)(t-s)\right)\,{\rm d}s~.$$
Conversely, it is easy to see that a solution of the previous 
integral equation is a solution of \eqref{eq-prop-W}.
To prove Proposition \ref{prop-W}, we set up a fixed point theorem for
contracting maps. We introduce the map $\Phi_{V}$, defined for $W$
bounded in $\Cc^0(\RR,\Qc_nX^{\sigma_0})$ by
$$\Phi_{V}(W)(t)=\int_0^\infty e^{\Qc_n A\Qc_n
s}\left((\Qc_nA\Pc_n)V(t-s)+\Qc_n F(V+W)(t-s)\right)\,{\rm d}s~.
$$

During the first part of the proof, we consider real valued functions, so that
the terms including the function $f$ are well defined. Let $\nu\in (0,1-\sigma_0]$ 
to be fixed later. Let $V=(v_1,v_2)$, we have $\Qc_n A\Pc_n V= (0,Q_n\gamma(x)P_n
v_2)$. Thus, using that $\gamma$ is of class $\Cc^1$,
$$\|\Qc_n A \Pc_n V\|_{X^{\nu+\sigma_0}} \leq C \|P_n
v_2\|_{H^{\nu+\sigma_0}} \leq C |\lambda_n|^{\nu/2}\|P_n
v_2\|_{H^{\sigma_0}}$$
where $(-\lambda_k)$ denotes the eigenvalues of the Laplacian operator $L$. 
Proposition \ref{prop-ana-1} and the bound on $V$ show that 
\begin{equation}\label{eq-demo-prop-W-1}
\left\| \int_0^\infty e^{\Qc_n A\Qc_n
s}(\Qc_nA\Pc_n) V(t-s)\,{\rm
d}s\right\|_{L^\infty(\RR,X^{\sigma_0})}~\leq~
C_1(\nu)|\lambda_n|^{\nu}\|V\|_{L^\infty(\RR,X^{\sigma_0})}
\end{equation}
as soon as there is $\varepsilon>0$ such that $1/t^{\frac{\nu\beta}{2}-\varepsilon}$ 
is integrable in a neighborhood of $+\infty$, that is for $\nu\beta>2$. 

For any function $Z(t)=(z_1,z_2)$ bounded in $X^{\sigma_0}$, we have
$\Qc_n F(Z)=(0,Q_n f(z_1))$. Since $H^{1+\sigma_0}(\Omega)$ is a normed 
algebra, using the polynomial growth of $f$ stated in \eqref{hyp-f-estim}, 
we get
$$\|\Qc_n F(Z)\|_{X^{\nu+\sigma_0}} = \|Q_n
f(z_1)\|_{H^{\nu+\sigma_0}}\leq \frac C {|\lambda_n|^{1-\nu}}
\|Q_n f(z_1)\|_{H^{1+\sigma_0}}\leq
\frac{C(1+\|Z\|_{X^{\sigma_0}}^p)}{|\lambda_n|^{1-\nu}}~.$$
Once again, Proposition \ref{prop-ana-1} shows that
\begin{equation}\label{eq-demo-prop-W-2}
\left\| \int_0^\infty e^{\Qc_n A\Qc_n
s}\Qc_n F(Z)(t-s)\,{\rm d}s\right\|_{ L^\infty(\RR,
X^{\sigma_0})}~\leq~\frac{C_2(\nu)(1+\|Z\|_{L^\infty(\RR,X^{\sigma_0}
) }^p ) } {|\lambda_n|^{ 1-\nu}}
\end{equation}
as soon as $\nu\beta>2$.

In the same way, using the control of $f'_u$ stated in \eqref{hyp-f-estim}, 
we prove that, if $\nu\beta>2$,
\begin{align}
&\left\| \int_0^\infty e^{\Qc_n A\Qc_n
s}\Qc_n \left(F(Z)-F(Z')\right)(t-s)\,{\rm
d}s\right\|_{L^\infty(\RR,X^{\sigma_0})}~\nonumber\\
&~~~~~~~~~~~~\leq~\frac{C_3(\nu)\left(1+\|Z\|_{
L^\infty(\RR, X^{\sigma_0})}^{p-1}+\|Z'\|_{
L^\infty( \RR,X^{\sigma_0})}^{p-1}\right)}{|\lambda_n|^{
1-\nu}}\|Z-Z'\|_{ L^\infty(\RR, X^{\sigma_0})}
~.\label{eq-demo-prop-W-3}
\end{align}
Gathering \eqref{eq-demo-prop-W-1}, \eqref{eq-demo-prop-W-2} and
\eqref{eq-demo-prop-W-3}, we obtain that, for any real functions $V$ and $W$
with 
$\|V\|_{ L^\infty( \RR,X^{\sigma_0})}\leq M$ and $\|W\|_{ L^\infty(\RR,
X^{\sigma_0})}\leq M'$, $\Phi$ is well defined, bounded and locally lipschitzian.
To apply the fixed point theorem for contracting maps, we need that the
Lipschitz constant is smaller than $1$, which is implied by
\begin{equation}\label{eq-demo-prop-W-4}
\frac{C_3(\nu)\left(1+M^{p-1}+{M'}^{p-1}\right)}{
|\lambda_n|^{ 1-\nu}}~\leq ~\frac 12
\end{equation}
and that $\Phi_V$ maps $B_{M',0}(\Qc_nX^{\sigma_0})$ into itself, 
which is implied by
\begin{equation}\label{eq-demo-prop-W-5}
C_1(\nu)|\lambda_n|^{\nu} M +
\frac{C_2(\nu)(1+M^p+{M'}^p)) } {|\lambda_n|^{ 1-\nu}} ~\leq~ \frac{M'}2~.
\end{equation}
To this end, we choose $\nu>0$ such that $(p-1)\nu < 1-\nu$ and we fix
$M'$ to be equal to $4C_1(\nu)|\lambda_n|^\nu M$, so that the first term 
of \eqref{eq-demo-prop-W-5} is smaller but satisfies the same growth 
than to the bound $M'/2$. Then, since $\lambda_n$
goes to $+\infty$ when $n\rightarrow +\infty$, one can find $n$ large enough
such that \eqref{eq-demo-prop-W-4} and \eqref{eq-demo-prop-W-5} hold, since 
$M=o(M')$ and ${M'}^{p-1}=o(|\lambda_n|^{1-\nu})$ when $n$ goes to 
$+\infty$. Taking $n$ larger if needed, the bounds $n\geq n_0$ and $M'\geq 
M'_0$ are easily fulfilled. 

It remains to check that we can find $\nu$ satisfying all the required conditions.
The bound $(p-1)\nu < 1-\nu$ is equivalent to $\nu<1/p$. It is compatible with $\nu\beta>2$ 
since we assumed $\frac{2}{\beta}<\frac{1}{p}$. Moreover, we also need that
$\nu\in (0,1-\sigma_0]$, which is possible since $\frac{2}{\beta}< 1-\sigma_0$.

\vspace{1mm}

Now, we extend our functions in a complex strip. By the previous bounds, if
it is real, $(v_1+w_1)(x,t)$ always stays smaller than $\kappa(M+M')$
where $\kappa$ is the is the injection constant $\|\cdot\|_{L^\infty}\leq \kappa
\|\cdot\|_{H^{1+\sigma_0}}$. Since $f$ is analytic, it has a holomorphic
extension in a complex neighborhood of the real interval
$[-\kappa(M+M');\kappa(M+M')]$. Thus, one can also consider functions $V$ and
$W$ with small imaginary parts in $X^{\sigma_0}$. All the above
estimates extend by continuity in this complex strip and, since
\eqref{eq-demo-prop-W-4} and \eqref{eq-demo-prop-W-5} contain some margin, for
$\delta,\delta'>0$ small enough, $\Phi_V$ can be extended as a contraction map
from $B_{M',\delta'}(\Qc_nX^{\sigma_0})$ into itself for all $V\in
B_{M+\delta,\delta}(\Pc_nX^{\sigma_0})$. Proposition \ref{prop-W} then
follows from the fact that $\Phi_V$ has a unique fixed point $W(V)$, which
corresponds to the unique solution of \eqref{eq-prop-W}.

To conclude, we notice that the above estimates also show that $\Phi_V(W)$ is
Lipschitz continuous with respect to $V$ and thus that the fixed point
$W(V)$ is Lipschitz continuous with respect to $V$. 
To obtain that the fixed point $W(V)$ depends holomorphically of $V$, we have to
show that the map $\Phi$ is holomorphic with respect to $V$ and $W$.
Then, one can conclude by using the implicit function theorem. To show that
$\Phi$ is holomorphic, we have to show that it has Fr\'echet derivatives. This
can be obtained by using the fact that $f$ is holomorphic and arguments similar
to the above ones.
\end{demo}

The proof of Proposition \ref{prop-W} also yields the property of {\it finite
determining modes}.
\begin{prop}\label{prop-W2}
Assume that Hypothesis (ii) of Theorem \ref{th-analytic} holds and let $M$ be 
given. Let $\sigma_0>0$ for $d=2$ or $\sigma_0>1/2$ for $d=3$. Then there 
exists $n\in\NN$ such that the following holds. Let $U_1(t)$ and $U_2(t)$ be two 
global solutions of \eqref{eq-analytic} such that $\|U_i(t)\|_{X^{\sigma_0}}\leq 
M$ for all $t\in\RR$. If $\Pc_nU_1(t)=\Pc_nU_2(t)$ for all times $t\in\RR$, then
$U_1(t)\equiv U_2(t)$ for all $t$.
\end{prop}
\begin{demo}
We consider the projections $V_i=\Pc_n U_i$ and $W_i=\Qc_n U_i$. For any $n$,
we have that $\|V_i(t)\|_{X^{\sigma_0}}\leq M$ and
$\|W_i(t)\|_{X^{\sigma_0}}\leq M$ for all $t\in\RR$. We argue as in the
proof of Proposition \ref{prop-W} with $M=M'_0$ and $n_0=0$. The fixed point
argument of Proposition \ref{prop-W} can be applied for $n$ large enough. It
shows that, if the low frequencies $V_1=V_2$ are known, then there is only one
possible high frequencies part $W(t)$, solution of \eqref{eq-prop-W}. Since
\eqref{eq-proj2} and \eqref{eq-proj} are equivalent, this unique function $W$
must be equal to $W_1$ as well as $W_2$. So $W_1=W_2$ and thus $U_1=U_2$.
\end{demo}

\vspace{3mm}

{\noindent \bf $\bullet$ Step 4: End of the proof of Theorem 
\ref{th-analytic}.}\\[1mm]
 Let $U(t)$ be
the mild solution of Theorem \ref{th-analytic} and let $U=V+W$ be the splitting
of \eqref{eq-proj2} for some $n\in\NN$. We have that $U(t)$ is uniformly bounded
in $X^{\sigma_0}$ by some constant $M$, thus, for all $n\in\NN$, $V(t)$
and $W(t)$ are also bounded by $M$, independent of $n$. From now on, we fix
$n$, $\delta$, $M'$ and $\delta'$ as prescribed by Proposition \ref{prop-W} for
such $M$, for $M'_0=M$ and $n_0=0$. We have the existence of a lipschitzian and
holomorphic map $\tilde V\mapsto W(\tilde V)$ defined in a neighborhood of $V$.
We consider the {\it ordinary differential equation} 
\begin{equation}\label{eq-demo-th-analytic}
\left\{\begin{array}{ll}
\partial_s \tilde V(s)=(\Pc_n A \Pc_n) \tilde V(s) + \Pc_n A\Qc_n W(\tilde V) +
\Pc_n F(\tilde V+W(\tilde V))\\
\tilde V(0)(t)=V(t) \in\Cc^0(\RR,\Pc_n X^{\sigma_0})
\end{array}\right.
\end{equation}
defined in the {\it Banach space} $\Cc^0(\RR,\Pc_n X^{\sigma_0})$ of finite dimension. Notice
that \eqref{eq-demo-th-analytic} is really an ODE since $\Pc_n$ has finite rank
and thus $(\Pc_n A)$ is a bounded operator. Proposition \ref{prop-W} shows
that $\tilde V\mapsto W(\tilde V)$ is lipschitzian and holomorphic in a
neighborhood of the initial data $V(t)$. It also ensures that $f$, and thus $F$,
are holomorphic in the complex set where $\tilde V+W(\tilde V)$ takes values. As
a consequence, by the classical theory of ODE's in Banach spaces,
\eqref{eq-demo-th-analytic} admits a unique solution $\tilde V(s)$ for
small $s\in \CC$, $|s|\leq \epsilon$, and this solution is holomorphic with
respect to $s$. Notice that the construction is just made such
that, for any $t\in\RR$ and for $s\in[-\epsilon,\epsilon]$ real,
$s\mapsto\tilde U(s)(t)=(\tilde V(s)+W(\tilde V(s)))(t)$ is a solution of
\eqref{eq-analytic} with $\Pc_n \tilde U(s)=\tilde V(s)$. By
uniqueness of the solution of \eqref{eq-demo-th-analytic}, using the translation
invariance, we have $\tilde U(0)(t+s)=\tilde U(s)(t)$. Thus $t\mapsto\tilde
U(0)(t)$ is a mild solution of \eqref{eq-analytic} with $\Pc_n (\tilde
U(0)(t))=V(t)$. Due to Proposition \ref{prop-W2}, $\tilde U(0)(t)=U(t)$ for all
$t\in\RR$. Since, for small $s\in[-\epsilon,\epsilon]$ real, we have
$U(t+s)=\tilde U(0)(t+s)=\tilde U(s)(t)$, we get that $s\mapsto U(t+s)$ is an
analytic function and thus $t\mapsto U(t)\in X^{\sigma_0}$ is an analytic
function.

\subsection{Proof of Theorem \ref{th1}}

Due to Proposition \ref{prop-UCP-implique-th}, we only need to obtain a unique 
continuation property. The difference with Theorem \ref{th2} is that the 
geometric background is quite general and we cannot use Theorems 
\ref{th-hormander-global} or \ref{th-hormander-global-2}. Our goal is thus to 
use the analytic unique continuation of \cite{Rob-Zui} stated in Theorem 
\ref{th-RZ}.

Let $W=(w,w_t)$ be a globally 
bounded solution of \eqref{eq} as in Propositions \ref{prop-asymp-compact-1} 
and \ref{prop-asymp-compact-2}. We have that $W(t)$ is globally bounded in 
$X^{\sigma_0}$ for $t\in\RR$ and $\sigma_0>0$. 
We want to apply Theorem \ref{th-analytic}. The assumption $\beta>2p\geq 2$ 
is already contained in the assumptions of Theorem \ref{th1}. Since we work here
with $\Omega$ of dimension $d=2$, we may choose $\sigma_0$ as small as wanted and in
particular we can fulfill the condition $\beta>2/(1-\sigma_0)$. If $\gamma$ is of class
$\Cc^1$, we may directly apply Theorem \ref{th-analytic}  and obtain 
that $W(t)$ is analytic in time. If $\gamma$ is not smooth, Assumption b) of Theorem \ref{th1}
provides a $\Cc^1$ damping $\tilde\gamma$ with a smaller support but with the same
decay properties for the corresponding damped wave semigroup. Since the energy of $W(t)$ is 
constant, $W(t)$ vanishes in the support of the damping $\gamma$
and satisfies
\begin{equation}\label{eq-demo-th1}
\partial_{tt}w(x,t)=\Delta w(x,t) - \alpha w(x,t) - f(x,w(x,t))~.
\end{equation}
We may thus replace $\gamma$ by the regular damping $\tilde\gamma$ and $W(t)$ is still
a solution of the corresponding damped (or free) wave equation. 
We apply Theorem \ref{th-analytic} in this setting and still get that $W(t)$ is analytic in time.

Since $w$ is smooth with respect to $t$ and $f$ is smooth, \eqref{eq-demo-th1} yields that 
$(\Delta - \alpha) w$ is in $L^2(\Omega)$ and thus $w$ is in $H^2(\Omega)$. 
We differentiate the above equation to obtain 
\begin{align*}
(\Delta-\alpha)^2 w&=(\Delta-\alpha)(\partial^2_{tt} w + f(x,w))\\
&=\partial^2_{tt}(\Delta-\alpha) w + (\Delta-\alpha) f(x,w)\\
&= \partial_{tttt}^4 w + f'_u(x,w)\partial^2_{tt}w
+ f''_{uu}(x,w)|\partial_t w|^2 + (\Delta-\alpha) f(x,w)
\end{align*}
showing that $w$ belongs to $H^4(\Omega)$. The process can be used as many 
times as wanted, showing that $w(x,t)$ is also smooth with respect to $x$.
Thus, the coefficients of \eqref{eq-UCP} are smooth in $x$ and analytic in $t$ 
and the unique continuation property of Theorem \ref{th-RZ} applies. Then 
Theorem \ref{th1} is a direct consequence of Proposition 
\ref{prop-UCP-implique-th}.

\section{Application 4: the disk with many holes}\label{section-appli4}

In Section \ref{section-appli3}, we have proved the semi-uniform stabilization 
for the semilinear damped wave equation in the case of the disk with two holes. 
In Figure \ref{fig-appli3}, we have seen that if the disk has three holes or 
more, Theorem \ref{th2} does 
not apply. In this case we assume that $u\mapsto f(x,u)$ is analytic and 
apply Theorem \ref{th1}. Notice that we may take $\gamma\in L^\infty$, since
the second part of Assumption b) of Theorem \ref{th1} is satisfied (see the geometric 
conditions in Theorem \ref{th-decay-disque-trous}).
Once again, we obtain an estimation of the decay which 
is better than the one given by Proposition \ref{prop1} because we follow the 
idea of the remark below Lemma \ref{lemma-appli2}. The details are left to the 
reader.
\begin{theorem}\label{th-appli4}
Let $\Oc\subset\RR^2$ be a smooth convex bounded open set. For $i=1\ldots p$, 
let $O_i\subset \Oc$ be smooth strictly convex obstacles satisfying:
\begin{enumerate}[(a)]
\item the obstacles are disjoint: $O_i\cap O_j=\emptyset$ for $i\neq j$,
\item no obstacle is in the convex hull of two others, that is that {\it 
convhull}$(O_i\cup O_j)\cap O_k =\emptyset$ for $i,j,k$ different,
\item if $\kappa$ denotes the infimum of the curvatures of the 
boundaries $\partial O_i$ of the obstacles and $L$ the minimal distance between 
two obstacles, and assume that $\kappa L>p$.  
\end{enumerate}
Let $\Omega=\Oc\setminus (\cup_i O_i)$ be the convex domain with 
holes and let $\alpha>0$. Assume moreover that 
\begin{enumerate}[(a)]
\item[(d)] the damping $\gamma\in L^\infty(\Omega,\RR_+)$ is strictly positive in a 
neighborhood of the exterior boundary $\partial \Oc$,
\item[(e)] the nonlinearity $f$ is smooth, satisfies \eqref{hyp-f-estim} and 
\eqref{hyp-f-sign} and $u\mapsto f(x,u)$ is analytic.
\end{enumerate}

Then, any solution $u$ of \eqref{eq} satisfies 
$$\|(u,\partial_t u)(t)\|_{H^1_0\times L^2}  ~\xrightarrow[~~t\longrightarrow 
+\infty~~]{}~0~.$$

Moreover, there exists $\tilde\lambda$ such that, for any $R$ and $\sigma\in 
(0,1]$, there exists $C_{R,\sigma}$ such that, for any solution $u$ with 
$U_0\in X^\sigma$, 
$$ \|(u_0,u_1)\|_{H^{1+\sigma}\times H^\sigma}\leq 
R~~\Longrightarrow ~~ \|(u,\partial_t u)(t)\|_{H^1_0\times L^2} \leq 
C_{R,\sigma}e^{-\sigma\tilde \lambda{t}^{1/3}}~.$$
\end{theorem}


\section{Application 5: Hyperbolic surfaces}\label{section-hyperbo}
In the case of hyperbolic surfaces, some recent results in Jin \cite{Lin} 
following the fractal uncertainty principle ideas of Bourgain-Dyatlov \cite{BourgainDyatlov} 
give a very good decay for any non trivial damping. In our nonlinear setting, the application 
of our previous results give the following result.
\begin{theorem}\label{th-appli5}
Let $\mathcal{M}$ be a compact connected hyperbolic surface
with constant negative curvature -1. Assume that 
\begin{enumerate}[(a)]
\item the damping $\gamma\in L^\infty(\Omega,\RR_+)$ is non zero and $\alpha>0$,
\item the nonlinearity $f$ is smooth, satisfies \eqref{hyp-f-estim} and 
\eqref{hyp-f-sign} and $u\mapsto f(x,u)$ is analytic.
\end{enumerate}

Then, any solution $u$ of \eqref{eq} satisfies 
$$\|(u,\partial_t u)(t)\|_{H^1_0\times L^2}  ~\xrightarrow[~~t\longrightarrow 
+\infty~~]{}~0~.$$

Moreover, there exists $\tilde\lambda$ such that, for any $R$ and $\sigma\in 
(0,1]$, there exists $C_{R,\sigma}$ such that, for any solution $u$ with 
$U_0\in X^\sigma$, 
$$ \|(u_0,u_1)\|_{H^{1+\sigma}\times H^\sigma}\leq 
R~~\Longrightarrow ~~ \|(u,\partial_t u)(t)\|_{H^1_0\times L^2} \leq 
C_{R,\sigma}e^{-\sigma\tilde \lambda\sqrt{t}}~.$$
\end{theorem}
The result follows from an application of Theorem \ref{th1}. 
The decay of the semigroup can be found in \cite{Lin}. Again, it is quite sure that 
we could avoid the loss and obtain a better decay $e^{-\sigma\tilde \lambda t}$ 
by following the arguments inside of the proof in \cite{Lin} for a slightly modified operator.

Note that the results in the references involve the case $\alpha=0$, 
but similar result can be obtained for $\alpha\geq 0$ for the linear semigroup 
as in Lemma \ref{lemma-appli2}. We also refer to other results with pressure conditions 
Schenck \cite{Schenckpressure} following ideas of Anantharaman \cite{Anant}.
We also want to stress that the result of \cite{Lin} follows several deep progress 
in the subject for ergodic flow and with various assumptions on the damping, but it 
would be impossible to make a complete bibliography. We refer to the bibliography in 
\cite{Lin} for instance, or the survey \cite{NonenmacherICM} for a history of resolvent 
estimates that can lead to such result of damping.


\section{A uniform bound for the $H^2\times H^1$-norm}\label{section-th-norm-H2}
This section is devoted to the proof of Theorem \ref{th-norm-H2}. We will 
assume in the whole section that the conclusions of Theorem \ref{th1} hold.

The following lemma is very general and does not depend on the geometric 
setting. 
\begin{lemma}\label{lemme-th-norm-H2}
Let $u(t)$ be a solution of the damped wave equation \eqref{eq} with 
$\gamma$ of class $\Cc^1$ and $f$ of class $\Cc^2$ satisfying 
\eqref{hyp-f-estim} and \eqref{hyp-f-sign}. Assume moreover that $\Omega$ is of 
dimension $d=2$ or of dimension $d=3$ and in this last case assume that $p < 3$. 
Then, if $U_0=(u_0,u_1)$ belongs to 
$D(A)=(H^2(\Omega)\times H^1_0(\Omega))\times H^1_0(\Omega)$, then 
$U(t)=(u,\partial_t u)(t)$ also belongs to $D(A)$. Moreover, there exists 
$\beta>0$ such that, for all $R>0$, there exists $C(R)>0$ such that,
$$\|(u_0,u_1)\|_{H^2\times H^1} \leq R~~~\Longrightarrow~~~\|(u,\partial_t 
u)(t)\|_{H^2\times H^1} \leq C(R) (1+t)^\beta~.$$
For $d=2$, the exponent $\beta$ is as close to $1$ as wanted.
\end{lemma}
\begin{demo}
The proof of this result is classical. Assume that $U_0=(u_0,u_1)$ belongs to 
$D(A)=(H^2(\Omega)\times H^1_0(\Omega))\times H^1_0(\Omega)$, we deal with the 
Cauchy problem by classical arguments since the linear semigroup $e^{At}$ is 
well defined on $D(A)$, $H^2(\Omega)\subset L^\infty(\Omega)$ and $f(x,0)=0$ 
due to \eqref{hyp-f-sign}, preserving the Dirichlet boundary conditions. Thus, 
the solution $U(t)$ is locally well defined in $D(A)$.

We recall that the damped wave equation admits the physical energy as a 
Lyapounov function
$$E(U)=\int_\Omega \frac 12(|\grad u|^2+\alpha |u|^2 + |\partial_t u|^2) + 
V(x,u)~.$$
As noticed in Section \ref{section-basic}, this energy is non-increasing 
in time. Using the Sobolev embeddings and Assumptions \eqref{hyp-f-estim} 
and \eqref{hyp-f-sign}, we obtain that $\|U(t)\|_{H^1\times L^2}$ is 
uniformly bounded if $U_0$ belongs to a bounded set of $H^1_0(\Omega)\times 
L^2(\Omega)$ and in particular is bounded by a constant $C_1(R)$ if 
$\|(u_0,u_1)\|_{H^2\times H^1} \leq R$.

We introduce a energy of higher order
$$\Fc(U)=\frac 12 \int_\Omega (|\Delta u|^2 + \alpha|\grad u|^2 + 
|\grad \partial_t u|^2) - \int_\Omega f(x,u)\Delta u~.$$
Notice that this energy is well defined for $U(t)\in D(A)$. To be more precise, 
we have 
$$\left|\int   f(x,u)\Delta u \right|  \leq \frac 14 \int |\Delta u|^2 + \int 
|f(x,u)|^2$$
and since $H^1(\Omega)\subset L^{2p}(\Omega)$, $\int |f(x,u)|^2$ is controlled 
by $\|u\|_{H^1}^{2p}$ and thus by $C_1(R)^{2p}$. In particular, there exists 
$C_2(R)$ such that
$$ \frac 14 \|U(t)\|_{H^2\times H^1}^2 - C_2(R) ~\leq~\Fc(U(t)~\leq ~ \frac 34 
\|U(t)\|_{H^2\times H^1}^2 + C_2(R)~.$$
We have
\begin{align}
\partial_t \Fc(U)&=\int \Delta u\Delta \partial_tu + \alpha \grad 
u\grad\partial_t u +  \grad\partial_tu.\grad\partial_{tt}^2 u - 
 f'_{u}(x,u)\partial_t u \Delta u -  f(x,u) \Delta\partial_t u \nonumber\\
&= \int (\Delta u - \alpha u - f(x,u) -\partial_{tt}u)\Delta\partial_t u
 - f'_{u}(x,u)\partial_t u \Delta u \nonumber\\
&= \int \gamma(x)\partial_t u \Delta\partial_t u 
- f'_{u}(x,u)\partial_t u \Delta u \nonumber\\
&= \int - \gamma(x) |\grad\partial_t u|^2 - \partial_t u \grad \gamma\grad  
\partial_t u - f'_{u}(x,u)\partial_t u \Delta u \nonumber\\
&\leq \int |\partial_t u| |\grad \gamma| |\grad  
\partial_t u| + \int |f'_{u}(x,u)||\partial_t u||\Delta u| 
\label{eq-preuve-lemme}
\end{align}
The first term is bounded by $\|\grad\gamma\|_{L^\infty}\|U\|_{H^1\times 
L^2}\|U\|_{H^2\times H^1}$ and so by $C_3(R)\|U\|_{H^2\times H^1}$. The second 
term of \eqref{eq-preuve-lemme} is 
bounded by $C_1(R) \|f'_u(\cdot,u)\|_{L^\infty} \|U\|_{H^2\times H^1}$. We 
bound $\|f'_u(\cdot,u)\|_{L^\infty}$ as follows for $d=3$:
\begin{align*}
\|f'_u(\cdot,u)\|_{L^\infty}&\leq (1+\|u\|_{L^\infty})^{p-1} \leq 
(1+\|u\|_{H^{3/2+\varepsilon}})^{p-1}\\
&\leq (1+\|u\|_{H^1}^{1/2-\varepsilon}\|u\|_{H^2}^{1/2+\varepsilon})^{p-1}
\end{align*}
where $\varepsilon\in (0,1/2]$ can be chosen small enough so that 
$\theta=(p-1)(1/2+\varepsilon)<1$ since $p<3$. Thus the second term of 
\eqref{eq-preuve-lemme}  is bounded by $C_4(R)(1+\|U\|_{H^2\times 
H^1}^{1+\theta})$. We finally obtain that
$$ \partial_t \Fc(U) \leq C_5(R) (C_6(R) + \Fc(U))^\eta$$
with $\eta=(1+\theta)/2<1$. This show that $\Fc(U)\leq C(R) (1+t)^\delta$ with 
$\delta=1/(1-\theta)$.

In the case $d=2$, the bound of the second term of 
\eqref{eq-preuve-lemme} is of the type 
$(1+\|u\|_{H^1}^{1-\varepsilon}\|u\|_{H^2}^{\varepsilon})^{p-1}$
with $\varepsilon\in (0,1]$ as small as needed. Thus the growth of 
$\Fc(U)$ is of type $(1+t)^\delta$ with $\delta$ as close to $2$ as wanted. 
Since $\Fc(U)$ is equivalent to $\|U\|_{H^2\times H^1}^2$, we obtain the 
polynomial growth of Lemma \ref{lemme-th-norm-H2} with $\beta$ as close to $1$ 
as wanted.
\end{demo}

{ \noindent \emph{\textbf{Proof of Theorem \ref{th-norm-H2}:}}}
By the previous lemma, we know that the $H^2\times H^1-$norm of $U(t)$ has a at 
most polynomial growth. By assumption, we know 
that the $H^1\times L^2-$norm of $U(t)$ goes to zero faster than any 
polynomial decay. By interpolation, this shows that the norm 
$\|U(t)\|_{H^{1+\varepsilon}\times H^\varepsilon}$ is bounded for 
$\varepsilon\in (0,1)$. Since $H^{1+\varepsilon}(\Omega)$ is an algebra, we 
have that $f(\cdot,u)$ is uniformly bounded in $H^{1+\varepsilon}(\Omega)$. 
Thus, for any initial data satisfying $\|U_0\|_{H^2\times H^1}\leq R$, and for 
any $\varepsilon \in (0,1/2)$, $F(U(t))=(0,f(\cdot,u(t)))$ is uniformly bounded 
in $(H^{2+\varepsilon}(\Omega)\cap H^1_0(\Omega)\times 
H^{1+\varepsilon}_0(\Omega)$ by a constant $C(R)$. 
We use a last time the formula of variation of the constant
$$U(t)=e^{At}U_0 + \int_0^t e^{At} F(U(t-s)) ds$$
and the weak decay estimates $\|e^{At}\|_{\Lc(H^{2+\varepsilon}(\Omega)\cap 
H^1_0(\Omega)\times H^{1+\varepsilon}_0(\Omega),D(A))}\leq h(t)$ with $h(t)$ 
integrable on $[0,+\infty)$ to obtain that $\|U(t)\|_{D(A)}$ is uniformly 
bounded. {\hfill$\square$\\}

{\noindent \bf Remark:} in the case of the open book, the decay of 
the semigroup is polynomial and Theorem \ref{th-norm-H2} does not apply. To 
adapt the above arguments and still get integrability where it is needed, we 
should assume that the vanishing order $\beta$ is small enough. However, this 
constraint seems not compatible with $\gamma$ of class $\Cc^1$. It could be 
possible to find sharper arguments but this is not the purpose of this paper. 

\appendix

\section{Estimates of the resolvent and decay of the semigroup}

The decay rate of a linear semigroup $e^{At}$ is closely related to the control 
of the resolvent $(A-i\mu)$ with $\mu\in\RR$, that is the resolvent along the 
imaginary axis. A famous result of \cite{Ge}, \cite{Pruss}
and \cite{Huang} is as follows. 
\begin{theorem}{\bf Gearhart-Pr\"uss-Huang}\label{th-Huang}\\
Let $e^{At}$ be a $\Cc^0-$semigroup in a Hilbert space $X$ and assume that
there exists a positive constant $M>0$ such that $\|e^{At}\|_{\Lc(X)}\leq M$ for
all $t\geq 0$. Then there exist $C>0$ and $\lambda >0$ such that
$$\forall U\in X~,~~\|e^{At}U\|_{X} \leq Ce^{-\lambda t} \|U\|_X$$
if and only if $i\RR\subset \rho(A)$ and 
$$\sup_{\mu\in\RR} \|(A-i\mu Id)^{-1}\|_{\Lc(X)}~<~+\infty~.$$
\end{theorem}

In the case of the weak stabilization, the resolvent $(A-i\mu 
Id)^{-1}$ is no more uniformly bounded for $\mu\in\RR$.
The rate of blow-up of this resolvent when $\mu\rightarrow \pm\infty$ is 
related to the decay of $e^{At}$ in $\Lc(D(A),X)$. A general relation has been 
obtained by Batty and Duyckaerts in \cite{Batty-Duyckaerts}. The first 
implication is the following.
\begin{theorem}\label{th-BD-1}
{\bf Batty-Duyckaerts (2008) \cite[Proposition 1.3]{Batty-Duyckaerts}.}\\
Let $e^{At}$ be a semigroup of operators on a space $X$ and assume that 
$$m(t)=\sup_{s\geq t} \|e^{At}(A-1)^{-1}\|_{\Lc(X)}$$
goes to $0$ as $t\rightarrow +\infty$. Then $i\RR$ belongs to the resolvent set 
of $A$ and there exist $\mu_0$ and $C>0$ such that 
$$\forall \mu\in \RR\text{ with }|\mu|\geq \mu_0~,~~
\|(A-i\mu)^{-1}\|_{\Lc(X)}\leq 1 + C m^{-1}_r \left(\frac 1{2(|\mu|+1)}\right)$$
where $m^{-1}_r$ is a right inverse of $m$, which maps $(0,m(0)]$ onto 
$[0,+\infty)$.
\end{theorem}
The second implication is more useful but is not optimal 
in general due to the logarithmic loss in $M_{\log}$ which is not expected.  
\begin{theorem}\label{th-BD}
{\bf Batty-Duyckaerts (2008) \cite[Theorem 1.5]{Batty-Duyckaerts}.}\\
Let $e^{At}$ be a semigroup of operators on a space $X$ such that
$\|e^{At}\|_{\Lc(X)}\leq C$ for all $t\geq 0$, and such that $i\RR\cap
\sigma(A)=\emptyset$. We set 
$$M(\mu)=\sup_{|\tau|\leq \mu} \|(A-i\mu)^{-1}\|_{\Lc(X)}$$
and 
$$M_{\log}(\mu)=M(\mu)[\ln(1+M(\mu))+\ln(1+|\mu|)]~.$$
Then, for any $k\in\NN\setminus\{0\}$, there exist $C_k$ and $T_k$, depending
only on $C$, $k$ and $M$, such that 
$$\forall t\geq T_k~,~~\left\|e^{At}(A-1)^{-k}\right\|\leq
\frac{C_k}{\left(M^{-1}_{\log}\left(\frac t{C_k}\right)\right)^k}$$
where $M^{-1}_{\log}$ is the inverse of $M_{\log}$ which maps 
$(M_{\log}(0),+\infty)$ onto $(0,+\infty)$.
\end{theorem}

If $X$ is a Hilbert space and $M$ polynomial, we can get rid of
the logarithmic term in $M_{\log}$ as proved in 
\cite{Bor-Tom} (see also \cite{Ana-Leau} in the framework of a damped wave 
system).
\begin{theorem}\label{th-BT}
{\bf Borichev-Tomilov (2010) \cite[Theorem 2.4]{Bor-Tom}.}\\
Let $e^{At}$ be a bounded $\Cc^0-$semigroup on a Hilbert space
$X$ with generator A such that $i\RR \cap \sigma(A)=\emptyset$. 
Then, for a fixed $\alpha>0$, the following conditions are equivalent:
\begin{enumerate}[(i)]
\item for large $\mu\in\RR$, $\|(A-i\mu)^{-1}\|_{\Lc(X)} = 
\Oc(|\mu|^{1/\alpha})$,
\item for large $t\geq 0$, $\|e^{At}A^{-1}\|_{\Lc(X)} = \Oc(
1/{t^{\alpha}})$,
\item for all $x\in H$ and for large $t\geq 0$, $\|e^{At}A^{-1}x\|_{X} = 
o(1/{t^{\alpha}})$.
\end{enumerate}
\end{theorem}

\section{Estimates of the resolvent of abstract damped wave 
equations}\label{section-appendix1}

In this section, we consider an abstract damped wave equation. Let $H$ be 
Hilbert spaces and let $L:D(L)\rightarrow H$ be a positive self-adjoint 
operator 
with compact resolvent. Let $B\in\Lc(H)$ be a damping operator which is 
bounded, self-adjoint and non-negative. We set $X=D(L^{1/2})\times H$ and 
$$A=\left(\begin{array}{cc} 0 & Id \\ -L & -B 
\end{array}\right)~~~~~~D(A)=D(L)\times D(L^{1/2})~.$$
For $\mu\in\RR$, we also introduce the operator $P_B(\mu):D(L)\rightarrow H$ 
defined by
$$P_B(\mu)=-L-i\mu B +\mu^2 Id~.$$
In Theorems \ref{th-BD} and \ref{th-BT}, we have seen the importance of 
estimating the resolvent $(A-i\mu)^{-1}$. In this section, we recall the 
equivalence with estimating $P_B(\mu)^{-1}$, which is often more convenient. 
This type of arguments is very classical and may be found in many articles 
dealing with the stabilization of damped wave equations. We present them here 
for sake of completeness and because we will need to generalize most of them in 
the next appendix.

We begin with the resolvent $(A-i\mu)^{-1}$ for fixed $\mu\in\RR$, that is that 
we consider the low frequencies.
\begin{prop}\label{prop-resolvent}
Let $\mu\in\RR$, the three following propositions are equivalent
\begin{enumerate}[(i)]
\item $(A-i\mu Id)$ is invertible in $\Lc(X)$,
\item $P_B(\mu)$ is invertible in $\Lc(H)$,
\item for any $u\neq 0$ solution of $Lu=\mu^2 u$, we have $\langle 
Bu|u\rangle\neq 0$.
\end{enumerate}
\end{prop}
\begin{demo}
Let $U=(u_1,u_2)$ and $V=(v_1,v_2)$ be vectors of $X=
D(L^{1/2})\times H$ such that $(A- i\mu)U=V$. We have equivalently
\begin{equation}\label{eq-demo-prop-resolvent}
\left\{\begin{array}{ll} 
u_2 -i\mu u_1 = v_1 & \text{ in }D(L^{1/2})\\
P_B(\mu) u_1 = v_2 + B v_1 + i \mu v_1&\text{ in }H~.
\end{array}\right. 
\end{equation}
Since $L$ has compact resolvent, the $P_B(\mu)$ is invertible if and only 
if its kernel is reduced to $\{0\}$ and if it does, $P_B(\mu)^{-1}$ is bounded 
from $H$ into $D(L)$. Thus \eqref{eq-demo-prop-resolvent} yields that 
(i)$\Leftrightarrow$(ii). Moreover, if $P_B(\mu)u=0$ with $u\neq 0$ then 
taking the scalar product with $u$ and considering the imaginary part, we get 
that $\langle 
Bu|u\rangle=0$ and so $B^{1/2}u=0$. Thus, $Bu=0$ and $Lu=\mu^2 u$ showing that 
(iii) fails. In the converse way, if we assume that (iii) fails, the 
corresponding solution $u$ also solves $P_B(\mu)u=0$ showing that (ii) fails.
\end{demo}

To study the high-frequencies, we have to estimate the behaviors for large 
$\mu$.
\begin{prop}\label{prop-estimP}
With the above notations, both estimations are equivalent 
\begin{enumerate}[(i)]
\item for large $\mu\in\RR$, $\|(A-i\mu)^{-1}\|_{\Lc(X)} = \Oc (M(|\mu|))$,
\item for large $\mu\in\RR$, $\|(P_B(\mu))^{-1}\|_{\Lc(H)}= \Oc( 
\frac{M(|\mu|)}{|\mu|})$.
\end{enumerate}
\end{prop}
\begin{demo}
The proof of (ii)$\Rightarrow$(i) is detailed in Proposition 
\ref{prop-2-appendix} below, adding projections on the high-frequencies. The 
implication stated here is simply the complete case $n=0$.

Let us show the converse implication.
Take $v_1=0$ in \eqref{eq-demo-prop-resolvent}. We have 
$u_2=i\mu u_1$ and $P_B(\mu)u_1 = v_2$. 
If (i) holds, that is $\|(A - i\mu )^{-1}\|_{\Lc(X)}\leq 
M(|\mu|)$, we must have in particular that $\|u_2\|_{H}\leq 
M(|\mu|)\|v_2\|_{H}$. Since $u_2=i\mu P_B(\mu)^{-1} v_2$, we obtain (ii).
\end{demo}

In some cases, to obtain an estimation of $\|(P_B(\mu))^{-1}\|_{\Lc(H)}$, it is 
more convenient to prove an observability estimate and to use the following 
proposition. Notice that this proposition yields a loss due to the term 
$f(\mu)^2$ in $M(\mu)$. This loss may sometimes be avoided but this may 
require an accurate study, based on particular dynamical properties of the 
geodesic flow.
\begin{prop}\label{prop-1-appendix}
We set
$$P(\mu)=-L+\mu^2 Id:=P_0(\mu)~.$$
Assume that there exist two positive functions $f$ and $g$ and $\mu_0\geq 0$ 
such that, for any $\mu$ with $|\mu|\geq \mu_0$ and 
any $u\in D(L)$, 
\begin{equation}\label{hyp-appendix}
\|u\|~\leq~ \frac{f(\mu)}\mu \|P(\mu)u\| + g(\mu) \|\sqrt{B} u\|~.
\end{equation}
Then, for any $\mu$ with $|\mu|\geq \mu_0$ and any 
$u\in D(L)$, 
\begin{equation}\label{result-prop-1-1}
\|u\|\leq \frac{M(\mu)}{|\mu|} \|P_{B}(\mu)u\|~,
\end{equation}
where
\begin{equation}\label{result-prop-1-3}
M(\mu)= 3 \max \left({f(\mu)}\,,\,f(\mu)^2\|\sqrt{B}\|_{\Lc(H)}^2\,,\, 
g(\mu)^2\right)~.
\end{equation}
\end{prop}
\begin{demo}
In fact, this proposition is simply Proposition 
\ref{prop-1-appendix-high-frequencies} below in the particular case $n=0$ 
that is $Q_n=Id$. We choose to copy this particular case in this section for 
clarity.   
\end{demo}

To finish, let us study the case where $L$ is replaced by $\tilde 
L=L+V$ where $V$ is a bounded non-negative operator, typically a potential or a 
linearized term. Using the previous propositions, we show that, if 
$M(\mu)=o(\mu)$, then the estimates for $L$ are equivalent to the estimates for 
$\tilde L$. Using Theorems \ref{th-BD} or \ref{th-BT}, we may obtain a relation 
between the decays of the semigroups.
\begin{prop}\label{prop-appendix-linearized}
We use the above notations and set 
$$\tilde A=\left(\begin{array}{cc} 0 & Id \\ -\tilde L & -B 
\end{array}\right)=\left(\begin{array}{cc} 0 & Id \\ -L-V & -B 
\end{array}\right)= A +\left(\begin{array}{cc} 0 & 0 \\ -V & 0 
\end{array}\right) ~.$$
Assume that for large $\mu\in\RR$, 
$\|(A-i\mu)^{-1}\|_{\Lc(X)} = \Oc (M(|\mu|))$ where $M(|\mu|)=o(\mu)$.
Then,  $\|(\tilde A-i\mu)^{-1}\|_{\Lc(X)} = \Oc (M(|\mu|))$ also holds for 
large $\mu$.
\end{prop}
\begin{demo}
If $\|(A-i\mu)^{-1}\|_{\Lc(X)} = \Oc (M(|\mu|))$ with $M(|\mu|)=o(\mu)$, then 
Proposition \ref{prop-estimP} shows that 
$\|P_B(\mu)^{-1}\|=\Oc(M(|\mu|)/\mu)=o(1)$. We set 
$$\tilde P_B(\mu)= -\tilde L-i\mu B +\mu^2 Id = P_B(\mu)-V~.$$
We have for large $\mu$
$$\tilde P_B(\mu)=P_B(\mu)(Id-P_B(\mu)^{-1} V)~,$$
showing that $\tilde P_B(\mu)$ is invertible for large $\mu$ since 
$P_B(\mu)^{-1}$ goes to $0$. In addition, it shows that the estimates for 
$\|\tilde P_B(\mu)^{-1}\|$ and $\|P_B(\mu)^{-1}\|$ are equivalent. Then the 
reverse implication of Proposition \ref{prop-estimP} finishes the proof. 
\end{demo}

\section{Estimates for the high-frequencies 
projections}\label{section-appendix2}
In Section \ref{section-reg-analytic}, we need to estimate the decay of the 
semigroup projected into the eigenspaces corresponding to the 
high frequencies of the Laplacian operator. This estimation is not direct 
in the cases where the projectors on high-frequencies do not commute with $A$. 
The purpose of this Section is to prove results yielding quickly to estimates 
of the decay of the high-frequencies by using the above results Theorem 
\ref{th-BD} and \ref{th-BT}. In particular, we generalized some results of 
Appendix \ref{section-appendix1} by showing that they hold uniformly with 
respect to cut-off frequency of the projection. Notice that this 
type of decay estimates for the high-frequency part of the solutions of the 
damped wave equation is related to Theorem 10 of \cite{Burq-Lebeau}, which shows 
that the eigenspaces corresponding to the high frequencies of the Laplacian 
operator are mainly preserved by the flow of the damped wave equation.

We use the notations of Appendix \ref{section-appendix1}.
Since $L$ is self-adjoint, positive and with compact resolvent, there exists an 
orthonormal basis 
$(\phi_k)_{k\geq 0}$ of eigenfunctions of $L$. We introduce the 
high-frequencies truncations $Q_n$, that are the
projectors on the space Span$\{\phi_k,k\geq n\}$
$$Q_nu=\sum_{k\geq n} \langle u | \phi_k\rangle \phi_k~.$$
We also introduce the sequence of high-frequencies projections 
$\Qc_n=(Q_n,Q_n)$ on $X$.

We consider in $Q_n H$ the operators
$$P_{Q_n B Q_n}(\mu)=Q_nP_B(\mu)Q_n= - L -i\mu Q_nBQ_n + \mu^2 Id~.$$
and the projection of $A$ on the high frequencies: that is, for any $V\in\Qc_n 
X$, 
$$\Qc_n A \Qc_n =\left( \begin{array}{cc} 0 & Id \\ -L &
-Q_n B Q_n\end{array} \right)~.$$
We prove a generalization of the classical implication of Proposition 
\ref{prop-estimP}, which is uniform with respect to the high-frequencies 
projections.

\begin{prop}\label{prop-2-appendix}
Assume that there exist a function $M(\mu)$,
uniformly positive, and $\mu_0\geq 0$ such that, for any $\mu$ with $|\mu|\geq 
\mu_0$, $n\in\NN$ and $u\in Q_nD(L)$, we have 
\begin{equation}\label{hyp-prop-2-appendix}
\|u\|_{H}\leq \frac{M(|\mu|)}{|\mu|} \|P_{Q_n B Q_n}(\mu)u\|_{H}~.
\end{equation}
Then, there exists $K>0$ such that, for all
$n\in\NN$ and all $|\mu|\geq \mu_0$,
$$\left\|\left(\Qc_n A\Qc_n - i\mu \right)^{-1}\right\|_{\Lc(\Qc_nX)}\leq
K M(|\mu|)~.$$
\end{prop}
\begin{demo}
Let $U=(u_1,u_2)$ and $V=(v_1,v_2)$ be vectors of $\Qc_n X=\Qc_n
(D(L^{1/2})\times H)$ such that $(\Qc_n A\Qc_n 
- i\mu)U=V$. We have 
\begin{equation}\label{eq-demo-prop-2-appendix-1}
\left\{\begin{array}{ll} 
u_2 -i\mu u_1 = v_1 & \text{ in }D(L^{1/2})\\
(-L+\mu^2 Id) u_1-i\mu Q_n B Q_n u_1 = v_2 + Q_n B v_1 
+ i \mu v_1&\text{ in }H
\end{array}\right. 
\end{equation}
We set $w=(P_{Q_n B Q_n}(\mu))^{-1}(v_1)$. Since $P_{Q_n B Q_n}(\mu)+ L 
+i\mu Q_n BQ_n=\mu^2 Id$, we have
\begin{align*}
w&=\frac 1{\mu^2}(P_{Q_n B Q_n}(\mu))^{-1}(P_{Q_n B Q_n}(\mu)v_1+ L v_1+i\mu
Q_n B Q_n v_1)\\
&= \frac 1{\mu^2} \left(v_1+(P_{Q_n B Q_n}(\mu))^{-1}(L v_1)+i\mu 
(P_{Q_n B Q_n}(\mu))^{-1} (Q_n B Q_n v_1)\right)
\end{align*}
and so
\begin{align}
\|w\|_{H}\leq & \frac 1{\mu^2}\|v_1\|_{H} + \frac 1{\mu^2} \|(P_{Q_n B
Q_n}(\mu))^{-1}\|_{\Lc(D(L^{-1/2}),H)} \|v_1\|_{D(L^{1/2})}\nonumber\\
&~~~+ \frac {1}\mu \|(P_{Q_n B
Q_n}(\mu))^{-1}\|_{\Lc(H)} \| B \|_{\Lc(H)} 
\|v_1\|_{H}~.\label{eq-demo-prop-2-appendix-2}
\end{align}
Let us estimate $\|(P_{Q_n B Q_n}(h))^{-1}\|_{\Lc(H,D(L^{1/2}))}$. We have
\begin{align*}
\|u\|_{D(L^{1/2})}^2&=\langle L u |u\rangle_{H}\\
&= \langle -P_{Q_n B Q_n}(\mu)u + \mu^2 u 
-i\mu Q_n B Q_n u|u\rangle_{H}\\
&\leq \Oc(\mu^2)\|u\|_{H}^2 + \|P_{Q_n B
Q_n}(\mu)u\|_{H}\|u\|_H\\
&\leq \left( \Oc(\mu^2)\|(P_{Q_n B
Q_n}(\mu))^{-1}\|^2_{\Lc(H,H)} + \|(P_{Q_n B
Q_n}(\mu))^{-1}\|_{\Lc(H,H)}\right) 
\|P_{Q_n B Q_n}(\mu)u\|_{H}^2 
\end{align*}
where the above estimation are independent of $n$. 
The estimate $\||(P_{Q_n B Q_n}(\mu))^{-1}\||_{\Lc(H,H)}$ is given by
Hypothesis \eqref{hyp-prop-2-appendix}, yielding
$$ \|u\|_{D(L^{1/2})}^2 \leq \left(M(|\mu|)^2 + \frac{M(|\mu|)}{|\mu|}\right) 
\|P_{Q_n B Q_n}(\mu)u\|_{H}^2~.$$
Using that $M$ is uniformly positive, $M(\mu)/\mu=o(M(\mu)^2)$ and so
$$\|(P_{Q_n B Q_n}(\mu))^{-1}\|_{\Lc(H,D(L^{1/2}))} = \Oc(M(|\mu|))~.$$
Since $P_{Q_n B Q_n}(\mu)$ 
defined from $D(L^{-1/2})$ in $H$ is the adjoint of $P_{Q_n B Q_n}(-\mu)$ 
defined from $H$ in $D(L^{1/2})$, we also have 
$$\|(P_{Q_n B Q_n}(\mu))^{-1}\|_{\Lc(D(L^{-1/2}),H)}=\Oc(M(|\mu|))~.$$
Coming back to \eqref{eq-demo-prop-2-appendix-2}, we obtain that 
$$\|w\|_{H}\leq K\frac{M(|\mu|)}{\mu^2}\|v_1\|_{D(L^{1/2})}~.$$
Considering \eqref{eq-demo-prop-2-appendix-1}, we have that 
$$ P_{Q_n B Q_n}(\mu)(u_1-i\mu w)=v_2+Q_n B v_1$$
and thus, due to \eqref{hyp-prop-2-appendix}, that
$\|u_1-i\mu w\|_{H}\leq \Oc(M(|\mu|)/\mu) \|V\|_X$. Together with the 
above estimate for $w$, we obtain that $\|u_1\|_{H}\leq \Oc(M(|\mu|)/\mu) 
\|V\|_X$ and, using the first equation of \eqref{eq-demo-prop-2-appendix-1}, 
that $\|u_2\|_{L^2}\leq \Oc(M(|\mu|))\|V\|_X$.

It remains to estimate $\|u_1\|_{D(L^{1/2})}$. To this end, we take the scalar 
product of second line of \eqref{eq-demo-prop-2-appendix-1} with $u_1$ and 
consider the real part to obtain
$$\|u_1\|_{D(L^{1/2})}^2-\mu^2 \|u_1\|_{H}^2\leq
\Oc(\mu)\|V\|_X\|u_1\|_{H}$$
and thus, due to the above estimates, $\|u_1\|_{D(L^{1/2})}^2 \leq 
\Oc(M(|\mu|)^2) \|V\|^2_X$.
\end{demo}

To obtain estimates as \eqref{hyp-prop-2-appendix}, it is convenient to 
generalize Proposition \ref{prop-1-appendix} to the case where high-frequencies 
projections appear. In this way, we can use the classical observability 
estimate without projections to study the decay of the high-frequencies 
semigroup.
\begin{prop}\label{prop-1-appendix-high-frequencies}
We set
$$P(\mu)=-L+\mu^2 Id:=P_0(\mu)~.$$
Assume that there exist two positive functions $f$ and $g$ and $\mu_0\geq 0$ 
such that, for any $\mu$ with $|\mu|\geq \mu_0$ and 
any $u\in D(L)$, 
\begin{equation}\label{hyp-appendix-high-frequencies}
\|u\|_H~\leq~ \frac{f(\mu)}\mu \|P(\mu)u\|_H + g(\mu) \|\sqrt{B} u\|_H 
~.
\end{equation}
Then, for any $\mu$ with $|\mu|\geq \mu_0$, any $n\in\NN$ and any $u\in 
Q_nD(L)$, 
\begin{equation}\label{result-prop-1-1-high-frequencies}
\|u\|_H\leq \frac{M(\mu)}{|\mu|} \|P_{Q_n B Q_n}(\mu)u\|_H~,
\end{equation}
where
\begin{equation}\label{result-prop-1-3-appendix-high-frequencies}
M(\mu)= 3 \max \left({f(\mu)}\,,\,f(\mu)^2\|\sqrt{B}\|_{\Lc(H)}^2\,,\, 
g(\mu)^2\right)~.
\end{equation}
\end{prop}
\begin{demo}
Let $u\in Q_n D(L)$. We have 
$$\langle P_{Q_n B Q_n}(\mu) u | u \rangle = -\langle L u| u
\rangle + \mu^2 \langle u|u\rangle - i \mu \langle B Q_n u|Q_n u \rangle~.$$
In particular, the imaginary part of $\langle P_{Q_n B 
Q_n}(\mu)u|u\rangle$ is $\mu \langle B Q_n u|Q_n
u\rangle$ and 
\begin{equation}\label{demo-eq-1-high-frequencies}
\forall u\in Q_n H~,~~\langle B Q_n u|Q_n u\rangle=\|\sqrt{B} u\|_H^2 \leq \frac
1{|\mu|} \|P_{Q_n B Q_n}(\mu) u\|_H \|u\|_H~.
\end{equation}
We write
\begin{equation}\label{demo-eq-2-high-frequencies}
\|P(\mu)u\|_H=\|P_{Q_n B Q_n}(\mu)u + i\mu Q_n B Q_n u\|_H \leq \|P_{Q_n B
Q_n}(\mu)u\|_H+|\mu|\|Q_n B Q_n u\|_H~.
\end{equation}
Then, we compute, for any $u\in Q_n H$
\begin{align*}
\|Q_n B Q_n u\|_H^2
&\leq \|B Q_n u\|_H^2 \leq \|\sqrt{B}\|^2_{\Lc(H)} \|\sqrt{B}Q_nu\|_H^2\\
& \leq \|\sqrt{B}\|^2_{\Lc(H)} \langle B Q_n u |Q_n u\rangle 
\end{align*}
and using \eqref{demo-eq-1-high-frequencies}, we obtain
$$ \|Q_n B Q_n u\|_H^2 \leq \frac {\|\sqrt{B}\|^2_{\Lc(H)}} {|\mu|} 
\|P_{Q_n B Q_n}(\mu)u\|_H\|u\|_H~.$$
Combining this result with \eqref{hyp-appendix-high-frequencies}, 
\eqref{demo-eq-1-high-frequencies} and 
\eqref{demo-eq-2-high-frequencies}, we get
\begin{align*}
 \|u\|_H&~\leq~\frac{f(\mu)}{|\mu|} \|P_{Q_n B
Q_n}(\mu)u\|_H + \frac{f(\mu)}{\sqrt {|\mu|}} {\|\sqrt{B}\|_{\Lc(H)}} 
\|P_{Q_n B Q_n}(\mu)u\|_H^{1/2}\|u\|_H^{1/2}\\
&~~~~~~~~~~~~~~~~~~+\frac{g(\mu)}{\sqrt {|\mu|}} \|P_{Q_n B
Q_n}(\mu)u\|_H^{1/2}\|u\|_H^{1/2}~.
\end{align*}
We bound the sum on the right by three times the largest of the three terms. 
Depending on which one is the largest one, we get three different bounds, which 
can be gathered in \eqref{result-prop-1-1-high-frequencies} and 
\eqref{result-prop-1-3-appendix-high-frequencies}.
\end{demo}



\end{document}